\documentclass[12pt]{article}
\usepackage{mathdots}

\usepackage{times,mathptmx,amsthm,amsmath}
\usepackage{pstricks,pst-node,graphics}
\usepackage{epsf}
\usepackage{graphicx}

\newtheorem{theorem}{Theorem}[section]

\theoremstyle{definition}

\theoremstyle{remark}

\numberwithin{equation}{section}

\newcommand{\Z}{{\bf Z}}
\newcommand{\R}{{\bf R}}
\newcommand{\vu}{{\bf u}}
\newcommand{\vv}{{\bf v}}
\newcommand{\vw}{{\bf w}}
\newcommand{\vzero}{{\bf 0}}
\newcommand{\ve}{{\bf e}}
\newcommand{\tT}{\frac{2}{3}}
\newcommand{\fT}{\frac{5}{3}}
\newcommand{\sT}{\frac{7}{3}}
\newcommand{\Un}{\ 1 \ }
\newcommand{\jd}{\ddots}
\newcommand{\ji}{\iddots}

\begin{document}

\begin{center}
\bf The combinatorics of frieze patterns and Markoff numbers \rm \\
\end{center}
\bigskip
\begin{center}
James Propp \\
Department of Mathematics, \\
University of Massachusetts Lowell \\
James\_Propp@uml.edu.ignorethis
\end{center}

\bigskip

\small

\noindent
{\sc Abstract}: This article, based on joint work with
Gabriel Carroll, Andy Itsara, Ian Le, Gregg Musiker,
Gregory Price, Dylan Thurston, and Rui Viana,
presents a combinatorial model based on perfect matchings
that explains the symmetries of the numerical arrays
that Conway and Coxeter dubbed frieze patterns.
This matchings model is a combinatorial interpretation
of Fomin and Zelevinsky's cluster algebras of type $A$.
One can derive from the matchings model
an enumerative meaning for the Markoff numbers,
and prove that the associated Laurent polynomials 
have positive coefficients
as was conjectured (much more generally) by Fomin and Zelevinsky.
Most of this research was conducted under the auspices of REACH
(Research Experiences in Algebraic Combinatorics at Harvard).

\normalsize

\bigskip

\section{Introduction} \label{sec-intro}

This article is part of a recent burst of activity
relating to what Sergey Fomin and Andrei Zelevinsky
have dubbed the ``Laurent phenomenon''
(described in greater detail below),
This phenomenon has algebraic, topological, and combinatorial aspects,
and it the third of these aspects that is developed here. 
In particular, we show how two examples of rational recurrences ---
the two-dimensional frieze patterns of Conway and Coxeter,
and the tree of Markoff numbers ---
relate to one another and to the Laurent phenomenon. 

A {\it Laurent polynomial\/} in the variables $x,y,\dots$
is a rational function in $x,y,\dots$ that can be expressed as
a polynomial in the variables $x,x^{-1},y,y^{-1},\dots$;
for example, 
the function $f(x) = (x^2+1)/x = x + x^{-1}$ is a Laurent polynomial,
but the composition $f(f(x)) = (x^4+3x^2+1)/x(x^2+1)$ is not.
The preceding example shows 
that the set of Laurent polynomials in a single variable
is not closed under composition.
This failure of closure also holds in the multivariate setting;
for instance, if $f(x,y)$, $g(x,y)$ and $h(x,y)$
are Laurent polynomials in $x$ and $y$,
then we would not expect to find
that $f(g(x,y),h(x,y))$ is a Laurent polynomial as well.
Nonetheless, it has been discovered that,
in broad class of instances (embraced as yet by no general rule),
``fortuitous'' cancellations occur that cause Laurentness to be preserved.
This is the ``Laurent phenomenon''
discussed by Fomin and Zelevinsky~\cite{FZL}.

Furthermore, in many situations where the Laurent phenomenon holds,
there is a positivity phenomenon at work as well,
and all the coefficients of the Laurent polynomials turn out to be positive.
In these cases, the functions being composed 
are Laurent polynomials with positive coefficients;
that is, they are expressions
involving only addition, multiplication, and division.
It should be noted that subtraction-free expressions
do not have all the closure properties one might hope for,
as the example $(x^3+y^3)/(x+y)$ illustrates:
although the expression is subtraction-free,
its reduced form $x^2 - xy + y^2$ is not.

Fomin and Zelevinsky have shown that a large part
of the Laurentness phenomenon
fits in with their general theory of cluster algebras.
This article will treat one important special case 
of the Laurentness-plus-positivity phenomenon,
namely the case associated with cluster algebras of type $A$,
discussed in detail in~\cite{FZY}.
The purely combinatorial approach taken 
in sections~\ref{sec-frieze} and~\ref{sec-sideways} of this article
obscures the links with deeper issues
(notably the representation-theoretic questions
that motivated the invention of cluster algebras),
but it provides the quickest and most self-contained way 
to prove the Laurentness-plus-positivity assertion in this case
(Theorem~\ref{thm-frieze}).
The frieze patterns of Conway and Coxeter,
and their link with triangulations of polygons,
will play a fundamental role,
as will perfect matchings of graphs derived from these triangulations.
(For a different, more algebraic way 
of thinking about frieze patterns, see~\cite{CaCh}.
For an extension of the result of this article into a broader setting,
see~\cite{BM}.)

In sections~\ref{sec-tropical} and~\ref{sec-variant} of this article,
two variations on the theme of frieze patterns are considered.
One is the tropical analogue,
which has bearing on graph-metrics in trees.
The other variant is based on a recurrence
that looks very similar to the frieze relation;
the variant recurrence
appears to give rise to tables of positive integers
possessing the same glide-reflection symmetry as frieze patterns,
but positivity and integrality are still unproved.

In section~\ref{sec-snake},
the graphs constructed in section~\ref{sec-frieze}
are viewed from a number of different perspectives
that relate them to existing literature.

In section~\ref{sec-markoff}, 
the constructions of sections~\ref{sec-frieze} and~\ref{sec-sideways} 
are specialized to a case
in which the triangulated polygons come from pairs of mutually visible points
in a dissection of the plane into equilateral triangles.
In this case, counting the matchings of the derived graphs
gives us an enumerative interpretation of Markoff numbers
(numbers satisfying the ternary cubic $x^2+y^2+z^2=3xyz$).
This yields a combinatorial proof of a Laurentness assertion
proved by Fomin and Zelevinsky in~\cite{FZL}
(namely a special case of their Theorem 1.10)
that falls outside of the framework of cluster algebras in the strict sense.
Fomin and Zelevinsky proved Theorem 1.10
by use of their versatile ``Caterpillar Lemma'',
but this proof did not settle the issue of positivity.
The combinatorial approach adopted here shows 
that all of the Laurent polynomials that occur in 
the three-variable rational-function analogue of the Markoff numbers ---
the ``Markoff polynomials'' --- 
are in fact positive (Theorem~\ref{thm-triples}).

Section~\ref{sec-other} concludes with some problems suggested 
by the main result of section~\ref{sec-markoff}.
One can try to generalize the combinatorial picture
by taking other dissections of the plane into triangles,
or one can try to generalize by considering other Diophantine equations.
There are tantalizing hints of a link between
the two proposed avenues of generalization,
but its nature is still obscure.

This work was supported by funds
from the National Science Foundation
and the National Security Agency.
It would have been impossible without
the deep insight and hard work 
of Gabriel Carroll, Andy Itsara, Ian Le, 
Gregg Musiker, Gregory Price, and Rui Viana
(all of whom were undergraduates
at the time the work was done)
and conversations with Dylan Thurston,
as well as earlier conversations with
Rick Kenyon and David Wilson.

\section{Triangulations and frieze patterns} \label{sec-frieze}

A {\it frieze pattern}~\cite{CoCo}
is an infinite array such as
$$
\begin{array}{cccccccccccccccccccc}
...&\Un&   &\Un&   &\Un&   &\Un&   &\Un&   &\Un&   &\Un&   &\Un&   &... \\[1.2ex]
...&   &\fT&   & 2 &   &\Un&   & 5 &   &\tT&   & 3 &   &\fT&   & 2 &... \\[1.2ex]
...& 4 &   &\sT&   &\Un&   & 4 &   &\sT&   &\Un&   & 4 &   &\sT&   &... \\[1.2ex]
...&   & 5 &   &\tT&   & 3 &   &\fT&   & 2 &   &\Un&   & 5 &   &\tT&... \\[1.2ex]
...&\Un&   &\Un&   &\Un&   &\Un&   &\Un&   &\Un&   &\Un&   &\Un&   &...
\end{array}
$$
consisting of $n-1$ rows, each periodic with period $n$,
such that all entries in the top and bottom rows are equal to 1
and all entries in the intervening rows satisfy the relation
$$
\begin{array}{ccc}
  & A &   \\[1.6ex]
B &   & C \\[1.6ex]
  & D &  
\end{array} \ \ : \ \ AD = BC-1 \ \ .
$$
The rationale for the term ``frieze pattern'' is that
such an array automatically possesses glide-reflection symmetry
(as found in some decorative friezes):
for $1 \leq m \leq n-1$, the $(n-m)$th row
is the same as the $m$th row, shifted rightward
(or, equivalently, leftward) by $n/2$ positions.
Hence the relation $AD=BC-1$
will be referred to below as the ``frieze relation''
even though its relevance to friezes and their symmetries
is not apparent from the algebraic definition.

Frieze patterns arose from Coxeter's study
of metric properties of polytopes,
and served as useful scaffolding for various sorts of metric data;
see e.g.~\cite{C1} (page 160),~\cite{C2}, and~\cite{C3}.
Typically some of the entries in a metric frieze pattern are irrational.
Conway and Coxeter completely classify
the frieze patterns whose entries are positive integers,
and show that these frieze patterns constitute 
a manifestation of the Catalan numbers. 
Specifically, there is a natural association
between positive integer frieze patterns
and triangulations of regular polygons with labeled vertices.
(In addition to~\cite{CoCo}, see the shorter discussion 
on pp.\ 74--76 and 96--97 of~\cite{CG}.)
Note that for each fixed $n$, any convex $n$-gon would serve here
just as well as the regular $n$-gon,
since we are only viewing triangulations combinatorially.

$$
\pspicture(0,-.75)(9,5.796)
\psset{xunit=1.5,yunit=1.5}
\pspolygon(2,0)(4,0)(5,1.732)(4,3.464)(2,3.464)(1,1.732)
\psline(2,0)(2,3.464)
\psline(2,3.464)(4,0)
\psline(4,0)(4,3.464)
\rput(1.8,-.4){6}
\rput(4.2,-.4){5}
\rput(5.4,1.732){4}
\rput(4.2,3.864){3}
\rput(1.8,3.864){2}
\rput(0.6,1.732){1}
\endpspicture
$$
\begin{center}
Figure 1. A triangulated 6-gon.
\end{center}

From every triangulation $T$ of a regular $n$-gon
with vertices cyclically labeled 1 through $n$,
Conway and Coxeter build an $(n-1)$-rowed frieze pattern
determined by the numbers $a_1,a_2,\dots,a_n$,
where $a_k$ is the number of triangles in $T$ incident with vertex $k$.
Specifically:
(1) the top row of the array is $\dots,1,1,1,\dots$;
(2) the second row (offset from the first) 
is $\dots,a_1,a_2,\dots,a_n,a_1,\dots$ (with period $n$); and
(3) each succeeding row (offset from the one before)
is determined by the frieze recurrence
$$
\begin{array}{ccc}
  & A &   \\[1.6ex]
B &   & C \\[1.6ex]
  & D &  
\end{array} \ \ : \ \ D = (BC-1)/A \ \ .
$$
E.g., the triangulation shown in Figure 1
determines the data $(a_1,\dots,a_6)=(1,3,2,1,3,2)$
and the 5-row frieze pattern
$$
\begin{array}{ccccccccccccccc}
...&1& &1& &1& &1& &1& &1& &1&... \\[1.1ex]
...& &1& &3& &2& &1& &3& &2& &... \\[1.1ex]
...&1& &2& &5& &1& &2& &5& &1&... \\[1.1ex]
...& &1& &3& &2& &1& &3& &2& &... \\[1.1ex]
...&1& &1& &1& &1& &1& &1& &1&...
\end{array}
$$

Conway and Coxeter show that the frieze relation,
applied to the initial rows $\dots,1,1,1,\dots$
and $\dots,a_1,a_2,\dots,a_n,\dots$,
yields a frieze pattern.
Note that implicit in this assertion
is the proposition that every entry in rows 1 through $n-3$ is non-zero
(so that the recurrence $D = (BC-1)/A$ never involves division by 0).
It is not a priori obvious that 
each of the entries in the array is positive
(since the recurrence involves subtraction)
or that each of the entries is an integer 
(since the recurrence involves division).
Nor is it immediately clear why
for $1 \leq m \leq n-1$, the $(n-m)$th row of the table
given by repeated application of the recurrence
should be the same as the $m$th row, shifted,
so that in particular the $n-1$st row, like the first row,
consists entirely of 1's.

These and many other properties of frieze patterns
are explained by a combinatorial model of frieze patterns
discovered by Carroll and Price~\cite{CP} 
(based on earlier work of Itsara, Le, Musiker, Price, and Viana;
see~\cite{I} and~\cite{M}, 
as well as~\cite{CILP}).
Given a triangulation $T$ as above, define a bipartite graph $G = G(T)$
whose $n$ black vertices $v$ correspond to the vertices of $T$, 
whose $n-2$ white vertices $w$ correspond to the triangular faces of $T$, 
and whose edges correspond to all incidences between vertices and faces in $T$
(that is, $v$ and $w$ are joined by an edge precisely if
$v$ is one of the three vertices of the triangle in $T$ associated with $w$).
For $i \neq j$ in the range $1,...,n$,
let $G_{i,j}$ be the graph obtained from $G$ 
by removing black vertices $i$ and $j$ and all edges incident with them,
and let $m_{i,j}$ be the number of perfect matchings of $G_{i,j}$
(that is, the number of ways to pair all $n-2$ of the black vertices
with the $n-2$ white vertices,
so that every vertex is paired to a vertex of the opposite color
adjacent to it).
For instance, for the triangulation $T$ of the 6-gon 
defined in Figure 1, the graph $G_{1,4}$ is as shown in Figure 2,
and we have $m_{1,4} = 5$ since this graph has 5 perfect matchings.
\bigskip
$$
\pspicture(1.5,0)(7.5,5.196)
\psset{xunit=.5,yunit=.5}
\psline(5,5.196)(6,10.392)
\psline(5,5.196)(6,0)
\psline(8,3.464)(6,10.392)
\psline(8,3.464)(12,0)
\psline(8,3.464)(6,0)
\psline(10,6.928)(6,10.392)
\psline(10,6.928)(12,10.392)
\psline(10,6.928)(12,0)
\psline(13,5.196)(12,10.392)
\psline(13,5.196)(12,0)
\pscircle[fillstyle=solid,fillcolor=black](6,0){.2}
\pscircle[fillstyle=solid,fillcolor=white](8,3.464){.2}
\pscircle[fillstyle=solid,fillcolor=black](12,0){.2}
\pscircle[fillstyle=solid,fillcolor=white](13,5.196){.2}
\pscircle[fillstyle=solid,fillcolor=black](12,10.392){.2}
\pscircle[fillstyle=solid,fillcolor=white](10,6.928){.2}
\pscircle[fillstyle=solid,fillcolor=black](6,10.392){.2}
\pscircle[fillstyle=solid,fillcolor=white](5,5.196){.2}
\endpspicture
$$
\begin{center}
Figure 2. The graph $G_{1,4}$.
\end{center}

\bigskip

\begin{theorem}[Carroll and Price~\cite{CP}] \label{thm-cp}
The Conway-Coxeter frieze pattern of a triangulation $T$ is just the array
\bigskip
$$
\begin{array}{ccccccccc}
\dots & m_{1,2} &          & m_{2,3} &         & m_{3,4} &         & m_{4,5} & \dots \\[2.0ex]
\dots &         & m_{1,3}  &         & m_{2,4} &         & m_{3,5} &         & \dots \\[2.0ex]
\dots & m_{n,3} &          & m_{1,4} &         & m_{2,5} &         & m_{3,6} & \dots \\[2.0ex]
\dots &         & m_{n,4}  &         & m_{1,5} &         & m_{2,6} &         & \dots \\[2.0ex]
      &  \vdots &          &  \vdots &         &  \vdots &         &  \vdots &
\end{array}
$$
\bigskip
where here as hereafter we interpret all subscripts mod $n$.
\end{theorem}

\noindent
Note that this claim makes the glide-reflection symmetry of frieze patterns
a trivial consequence of the fact that $G_{i,j} = G_{j,i}$.

\begin{proof}
Here is a sketch of the main steps of the proof:

(1) $m_{i,i+1} = 1$: This holds because there is a tree structure on
the set of triangles in $T$ that induces a tree structure on the set
of white vertices of $G$.  If we examine the white vertices of $G$,
proceeding from outermost to innermost, we find that we have no
freedom in how to match them with black vertices, when we keep in
mind that every black vertex must be matched with a white vertex.
(In fact, the same reasoning shows that $m_{i,j} = 1$ whenever
the triangulation $T$ contains a diagonal connecting vertices 
$i$ and $j$.)

(2) $m_{i-1,i+1} = a_i$: The argument is similar, except now we
have some freedom in how the $i$th black vertex is matched: it
can be matched with any of the $a_i$ adjacent white vertices. 

(3) $m_{i,j} \, m_{i-1,j+1} = m_{i-1,j} \, m_{i,j+1} - 1$:
If we move the 1 to the left-hand side, we can use (1) to write 
the equation in the form
$$m_{i,j} \, m_{i-1,j+1} + m_{i-1,i} \, m_{j,j+1} = m_{i-1,j} \, m_{i,j+1}.$$
This relation is a direct consequence of a lemma due to Eric Kuo
(Theorem 2.5 in~\cite{Kuo}),
which is stated here for the reader's convenience: 

{\it Condensation lemma:}
If a bipartite planar graph $G$ has 2 more black vertices than white vertices,
and the black vertices $a,b,c,d$ lie in cyclic order on some face of $G$, then 
$$m(a,c) m(b,d) = m(a,b) m(c,d) + m(a,d) m(b,c),$$
where $m(x,y)$ denotes the number of perfect matchings
of the graph obtained from $G$
by deleting vertices $x$ and $y$ and all incident edges.

(1) and (2) tell us that Carroll and Price's theorem
applies to the first two rows of the frieze pattern,
and (3) tells us (by induction) that the theorem 
applies to all subsequent rows.
\end{proof}

It should be mentioned that Conway and Coxeter
give an alternative way of describing the entries in frieze patterns,
as determinants of tridiagonal matrices.
Note that $m_{i-1,i+1} = a_i$
which equals the determinant of the 1-by-1 matrix
whose sole entry is $a_i$,
while $m_{i-1,i+2} = a_{i} a_{i+1} - 1$ 
which equals the determinant of the 2-by-2 matrix
$$\left( \begin{array}{cc}
a_i &  1 \\
 1  & a_{i+1} \end{array} \right).$$
One can show by induction using Dodgson's determinant identity
(for a statement and a pretty proof of this identity see~\cite{Z2})
that $m_{i-1,i+k}$ equals the Euler continuant $[a_i,\dots,a_{i+k-1}]$,
that is, the determinant of the $k$-by-$k$ matrix
with entries $a_i,\dots,a_{i+k-1}$ down the diagonal,
1's in the two flanking diagonals, and 0's everywhere else.
This is true for any array satisfying the frieze relation
whose initial row consists of 1's, whether or not it is a frieze pattern.
Thus, any numerical array constructed via the frieze relation
from initial data consisting of a first row of 1's
and a second row of integers
will be an array of integers,
since entries in subsequent rows are
equal to determinants of integer matrices.
(One caveat is in order here:
although the table of tridiagonal determinants 
always satisfies the frieze relation,
it may not be possible to compute the table
using just the frieze relation,
since some of the expressions that arise
might be indeterminate fractions of the form $0/0$.)
However, for most choices of positive integers $a_1,\dots,a_n$,
the resulting table of integers
will not be an $(n-1)$-rowed frieze pattern,
because some entries lower down in the table will be negative (or vanish).
Indeed, Conway and Coxeter show that 
every $(n-1)$-rowed frieze pattern whose entries are positive integers
arises from a triangulated $n$-gon in the fashion described above.

\section{The sideways recurrence and its periodicity} \label{sec-sideways}

Recall that any $(n-1)$-rowed array of real numbers
that begins and ends with rows of 1's 
and satisfies the frieze relation in between,
with all rows having period $n$,
qualifies as a frieze pattern.

Note that if the vertices $1,\dots,n$ of an $n$-gon lie on a circle
and we let $d_{i,j}$ be the distance between points $i$ and $j$,
Ptolemy's theorem on the lengths of the sides and diagonals
of an inscriptible quadrilateral gives us the three-term quadratic relation
$$d_{i,j} \, d_{i-1,j+1} + d_{i-1,i} \, d_{j-1,j} = d_{i-1,j} \, d_{i,j+1}$$
(with all subscripts interpreted mod $n$).
Hence the numbers $d_{i,j}$ with $i \neq j$,
arranged just as the numbers $m_{i,j}$ were,
form an $(n-1)$-rowed array that almost qualifies as a frieze pattern
(the array satisfies the frieze relation
and has glide-reflection symmetry 
because $d_{i,j}=d_{j,i}$ for all $i,j$,
but the top and bottom rows do not in general consist of 1's).
The nicest case occurs when the $n$-gon is a regular $n$-gon of side-length 1;
then we get a genuine frieze pattern
and each row of the frieze pattern is constant.

Another source of frieze patterns is an old result from spherical geometry:
the {\it pentagramma mirificum\/} of Napier and Gauss
embodies the assertion that the arc-lengths of the sides
in a right-angled spherical pentagram
can be arranged to form the middle two rows of a four-rowed frieze pattern.

Conway and Coxeter show that frieze patterns 
are easy to construct if one proceeds not from top to bottom
(since one is unlikely to choose numbers $a_1,\dots,a_n$ in the second row
that will yield all 1's in the $(n-1)$st row) but from left to right,
starting with a zig-zag of entries connecting the top and bottom rows
(where the zig-zag path need not alternate 
between leftward steps and rightward steps
but may consist of any pattern of leftward steps and rightward steps),
using the sideways frieze recurrence
$$
\begin{array}{ccc}
  & A &   \\[1.5ex]
B &   & C \\[1.5ex]
  & D &  
\end{array} \ \ : \ \ C = (AD+1)/B
$$
Although a priori one might imagine that repeating this recurrence
would lead one to non-integer rational numbers
whose numerators and denominators would get increasingly large
as one goes from left to right,
it turns out that the resulting pattern 
necessarily repeats with period $n$,
and that all the numbers that appear are whole numbers
(provided that all the entries in the initial zig-zag
are equal to 1).

E.g., consider the partial frieze pattern
$$
\begin{array}{ccccccccccc}
... & 1 &   & 1 &    & 1 &   & 1 &   & 1 & ... \\[1.5ex]
    &   & x &   & x' &   &   &   &   &   &     \\[1.5ex]
    & y &   & y' &   &   &   &   &   &   &     \\[1.5ex]
    &   & z &   & z' &   &   &   &   &   &     \\[1.5ex]
... & 1 &   & 1 &    & 1 &   & 1 &   & 1 & ... \end{array}
$$
Given non-zero values of $x$, $y$, and $z$, one can successively compute 
$y' = (xz+1)/y$, 
$x' = (y'+1)/x$, 
and 
$z' = (y'+1)/z$, 
obtaining a new zig-zag of entries $x',y',z'$
connecting the top and bottom rows.
It is clear that for generic choices of non-zero $x,y,z$, 
one has $x',y',z'$ non-zero as well,
so the procedure can be repeated, yielding further zig-zags of entries.
After six iterations of the procedure
one recovers the original numbers $x,y,z$
six places to the right of their original position
(unless one has unluckily chosen $x,y,z$
so as to cause one to encounter an indeterminate expression
of the form $0/0$ from the recurrence),
and if we specialize to $x=y=z=1$,
we get the 5-row frieze pattern associated with Figure 1.

To dodge the issue of indeterminate expressions of the form $0/0$,
we embrace indeterminacy of another sort by regarding $x,y,z$
as formal quantities, not specific numbers,
so that $x',y',z'$, etc.\ become rational functions of $x$, $y$, and $z$.
Then our recurrence ceases to be problematic.
Indeed, one finds that the rational functions that arise
are of a special kind,
namely, Laurent polynomials with positive coefficients.

We can see why this is so by incorporating weighted edges
into our matchings model.
Returning to the triangulated hexagon from section~\ref{sec-frieze},
associate the values $x$, $y$, and $z$
with the diagonals joining vertices 2 and 6,
vertices 2 and 5, and vertices 3 and 5, respectively.
Call these the formal weights of the diagonals.
Also assign weight 1 to each of the 6 sides of the hexagon;
see Figure 3.

$$
\pspicture(0,-.75)(9,5.796)
\psset{xunit=1.5,yunit=1.5}
\pspolygon(2,0)(4,0)(5,1.732)(4,3.464)(2,3.464)(1,1.732)
\psline(2,0)(2,3.464)
\psline(2,3.464)(4,0)
\psline(4,0)(4,3.464)
\rput(1.4,.666){1}
\rput(1.4,2.798){1}
\rput(3.0,3.664){1}
\rput(4.6,2.798){1}
\rput(4.6,.666){1}
\rput(3.0,-.2){1}
\rput(1.8,1.732){$x$}
\rput(2.8,1.732){$y$}
\rput(4.2,1.732){$z$}
\endpspicture
$$
\begin{center}
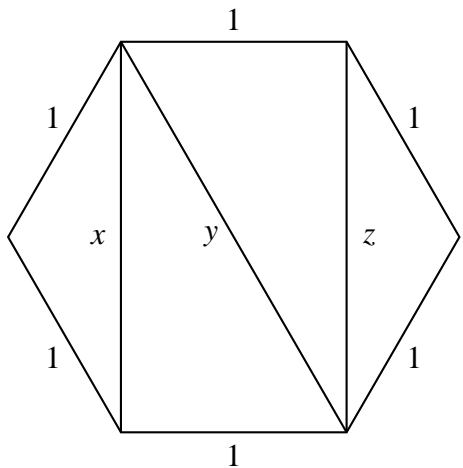
Figure 3. A triangulated 6-gon with edge-weights.
\end{center}

\bigskip

\noindent
Now construct the graph $G$ from the triangulation as before,
this time assigning weights to all the edges.
Specifically, if $v$ is a black vertex of $G$
that corresponds to a vertex of the $n$-gon
and $w$ is a white vertex of $G$
that corresponds to a triangle in the triangulation $T$
that has $v$ as one of its three vertices
(and has $v'$ and $v''$ as the other two vertices),
then the edge in $G$ that joins $v$ and $w$ should be assigned the weight
of the side or diagonal in $T$ that joins $v'$ and $v''$;
see Figure 4.

$$
\pspicture(1.5,0)(7.5,5.196)
\psset{xunit=.5,yunit=.5}
\psline(5,5.196)(6,10.392)
\rput(5.1,7.796){$1$} 
\psline(5,5.196)(6,0)
\rput(5.1,2.598){$1$} 
\psline(8,3.464)(6,10.392)
\rput(7.4,6.928){$1$} 
\psline(8,3.464)(12,0)
\rput(9.8,1.532){$x$} 
\psline(8,3.464)(6,0)
\rput(7.4,1.732){$y$} 
\psline(10,6.928)(6,10.392)
\rput(8,9.080){$z$} 
\psline(10,6.928)(12,10.392)
\rput(10.6,8.660){$y$} 
\psline(10,6.928)(12,0)
\rput(10.6,3.464){$1$} 
\psline(13,5.196)(12,10.392)
\rput(12.9,7.794){$1$} 
\psline(13,5.196)(12,0)
\rput(12.9,2.598){$1$} 
\psline(3,5.196)(5,5.196)
\rput(4,5.596){$x$} 
\psline(15,5.196)(13,5.196)
\rput(14,5.596){$z$} 
\pscircle[fillstyle=solid,fillcolor=black](6,0){.2}
\pscircle[fillstyle=solid,fillcolor=white](8,3.464){.2}
\pscircle[fillstyle=solid,fillcolor=black](12,0){.2}
\pscircle[fillstyle=solid,fillcolor=white](13,5.196){.2}
\pscircle[fillstyle=solid,fillcolor=black](12,10.392){.2}
\pscircle[fillstyle=solid,fillcolor=white](10,6.928){.2}
\pscircle[fillstyle=solid,fillcolor=black](6,10.392){.2}
\pscircle[fillstyle=solid,fillcolor=white](5,5.196){.2}
\pscircle[fillstyle=solid,fillcolor=black](3,5.196){.2}
\pscircle[fillstyle=solid,fillcolor=black](15,5.196){.2}
\endpspicture
$$
\begin{center}
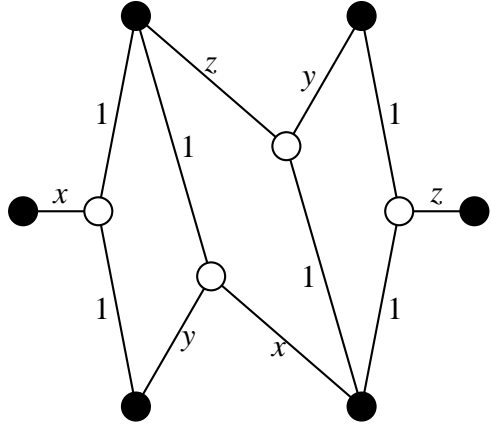
Figure 4. The graph $G_{1,4}$ with edge-weights.
\end{center}

\bigskip

\noindent
We now define $W_{i,j}$ as the sum of the weights of
all the perfect matchings of the graph $G_{i,j}$
obtained from $G$ by deleting vertices $i$ and $j$
(and all their incident edges),
where the weight of a perfect matching
is the product of the weights of its constituent edges,
and we define $M_{i,j}$ as $W_{i,j}$ divided by
the product of the weights of all the diagonals
(this product is $xyz$ in our running example);
e.g., $W_{1,4} = 1 + 2y + y^2 + xz$
and $M_{1,4} = x^{-1} y^{-1} z^{-1} + 2 x^{-1} z^{-1} + x^{-1} y z^{-1} + y^{-1}$.
These $M_{i,j}$'s, which are rational functions of $x$, $y$, and $z$,
generalize the numbers denoted by $m_{i,j}$ earlier,
since we recover the $m_{i,j}$'s from the $M_{i,j}$'s
by setting $x=y=z=1$.
It is clear that each $W_{i,j}$ is a polynomial with positive coefficients,
so each $M_{i,j}$ is a Laurent polynomial with positive coefficients.
And, because of the normalization (division by $xyz$),
we have gotten each $M_{i,i+1}$ to equal 1.
So the table of rational functions $M_{i,j}$
is exactly what we get by running our recurrence from left to right.
When we pass from $x,y,z$ to $x',y',z'$,
we are effectively rotating our triangulation by one-sixth of a full turn;
six iterations bring us back to where we started.

It is not hard to see that
the same approach works for any triangulation of an $n$-gon for any $n$,
and in this way we can prove:

\begin{theorem} \label{thm-frieze}
Given any assignment of formal weights
to $n-3$ entries in an $(n-1)$-rowed table
that form a zig-zag joining the top row (consisting of all 1's)
to the bottom row (consisting of all 1's),
there is a unique assignment of rational functions
to all the entries in the table so that the frieze relation is satisfied.
These rational functions of the original $n-3$ variables
have glide-reflection symmetry that gives each row period $n$.
Furthermore, each of the rational functions in the table
is a Laurent polynomial with positive coefficients.
\end{theorem}

Note that a zig-zag joining the top row to the bottom row
corresponds to a triangulation $T$ whose dual tree is just a path.
Not every triangulation is of this kind
(for instance, consider the hexagon shown in Figure 1
triangulated by diagonals joining vertices 2 and 4, 
vertices 4 and 6, and vertices 2 and 6).
In general, the entries in a frieze pattern
that correspond to the diagonals of a triangulation $T$
do not form a zig-zag path,
so it is not clear from the frieze pattern
how to extend the known entries to the unknown entries
(e.g., for the triangulation described in the parenthetical remark
in the previous sentence,
if one assigns respective weights $x$, $y$, and $z$
to the specified diagonals,
one obtains the partial frieze pattern
$$
\begin{array}{cccccccccccccccccccc}
...& 1 &   & 1 &   & 1 &   & 1 &   & 1 &   & 1 &   &... \\[1.2ex]
...&   & ? &   & x &   & ? &   & y &   & ? &   & z &... \\[1.2ex]
...& ? &   & ? &   & ? &   & ? &   & ? &   & ? &   &... \\[1.2ex]
...&   & y &   & ? &   & z &   & ? &   & x &   & ? &... \\[1.2ex]
...& 1 &   & 1 &   & 1 &   & 1 &   & 1 &   & 1 &   &...
\end{array}
$$
where the question marks refer to entries
whose values do not follow immediately
from the frieze recurrence).
In such a case, it is best to refer directly to the triangulation itself,
and to use a generalization of the frieze relation,
namely the (formal) Ptolemy relation~\cite{CP}
$$M_{i,j} \, M_{k,l} + M_{j,k} \, M_{i,l} = M_{i,k} \, M_{j,l}$$
where $i,j,k,l$ are four vertices of the $n$-gon listed in cyclic order.
(Conway and Coxeter~\cite{CoCo} 
give spatially extended versions of the frieze relation
that are equivalent to special cases of the Ptolemy relation.)
Since every triangulation of a convex $n$-gon 
can be obtained from every other by means of flips 
that replace one diagonal of a quadrilateral by the other diagonal
(an observation that goes back at least as far as 1936~\cite{W}),
we can iterate the Ptolemy relation so as to solve for all of the $M_{i,j}$'s
in terms of the ones whose values were given.

Up until now we have associated indeterminates
with the $n-3$ diagonals, but not the $n$ sides,
of our triangulated $n$-gon.
If we assign formal indeterminates to the sides as well as the diagonals
and carry out the construction of the edge-weighted graph $G_{i,j}$
(incorporating the $n$ new variables)
and the polynomial $W_{i,j}$
(the sum of the weights of all the perfect matchings of $G_{i,j}$),
and we define the Laurent polynomial $M_{i,j}$
as $W_{i,j}$ divided by the product of all $n-3$ diagonal-weights,
the proof of Theorem~\ref{thm-frieze} still goes through, 
and one sees that the $M_{i,j}$'s form an array
in which the top and bottom rows
contain the indeterminates associated with the sides of the $n$-gon
and the intervening rows satisfy the modified frieze relation
$$
\begin{array}{ccccccccl}
 X  &     &     &     &     &     &  Y  &   & \\
    & \jd &     &  A  &     & \ji &     &   & \\[1.2ex]
    &     &  B  &     &  C  &     &     & : & \ \ BC = AD + XY \\
    & \ji &     &  D  &     & \jd &     &   & \\
 Y  &     &     &     &     &     &  X  &   & \\
\end{array}
$$
where $X$ is the top entry in the diagonal containing $B$ and $D$ 
as well as the bottom entry in the diagonal containing $A$ and $C$
and $Y$ is the top entry in the diagonal containing $C$ and $D$
as well as the bottom entry in the diagonal containing $A$ and $B$. 
Each entry in this generalized frieze pattern
is a Laurent monomial in the $2n-3$ variables
in which the $n$ variables associated with sides of the $n$-gon
occur in only with non-negative exponents.

Our combinatorial construction
of Laurent polynomials associated with the sides and diagonals of an $n$-gon
is essentially nothing more than
the type $A$ case (more precisely, the $A_{n-3}$ case)
of the cluster algebra construction
of Fomin and Zelevinsky~\cite{FZY}.
The result that our matchings model yields,
stated in a self-contained way, is as follows:

\begin{theorem}
Given any assignment of formal weights
$x_{i,j}$ to the $2n-3$ edges of a triangulated convex $n$-gon,
(where $x_{i,j}$ is associated with the edge joining vertices $i$ and $j$),
there is a unique assignment of rational functions
to all $n(n-3)/2$ diagonals of the $n$-gon
such that the rational functions assigned 
to the four sides and two diagonals of any quadrilateral 
determined by four of the $n$ vertices satisfy the Ptolemy relation.
These rational functions of the original $2n-3$ variables
are Laurent polynomials with positive coefficients.
\end{theorem}

The formal weights are precisely the cluster variables
in the cluster algebra of type $A_{n-3}$,
and the triangulations are the clusters.
The periodicity phenomenon
is a special case of a more general periodicity conjectured by Zamolodchikov
and proved in the type $A$ case independently by Frenkel and Szenes
and by Gliozzi and Tateo; see~\cite{FZY} for details.


\section{Snake graphs} \label{sec-snake}

The bipartite graphs $G_{i,j}$ obtained in section~\ref{sec-frieze},
when shorn of their forced edges (edges that belong to every perfect matching)
and their forbidden edges (edges that belong to no perfect matching),
have a direct combinatorial construction as ``snakes'' of 4-cycles,
obtained by repeating the process of adding a new 4-cycle 
at the end of a snake.
More precisely: a snake of order 0 is just two vertices joined by an edge;
a snake of order 1 is a 4-cycle;
a snake of order 2 is a pair of 4-cycles sharing a single edge,
obtained by adjoining one 4-cycle to another along an edge;
and if one has a snake of order $k-1$
whose most recently added 4-cycle $C$ was adjoined along edge $e$,
one obtains a snake of order $k$
by adjoining a new 4-cycle
that shares some edge of $C$ other than $e$.
For example, Figure 5 shows a snake of order 6
obtained from a triangulated 9-gon
whose vertices are shown in black
(the two vertices represented by smaller black dots
are not part of the snake,
but they are included for clarity).
Given a triangulation $T$,
the only edges of $G_{i,j}$ that are neither forced nor forbidden
are those whose white endpoint
corresponds to a triangle in $T$
on the path of triangles joining vertices $i$ and $j$. 
These edges form a snake-graph
whose twists and turns mimic those of the path of triangles.

Up until now, we have used $n$ to denote the number of sides 
of the polygon being triangulated,
but in this section it will be more convenient 
to let $n+3$ denote the number of sides of the polygon,
and to make the additional assumption
that every triangle in the triangulation of the $(n+3)$-gon
occurs on the path of triangles from $i$ to $j$,
so that the graph $G_{i,j}$ has no forced or forbidden edges
and the snake-graph has order $n$.

$$
\pspicture(-2.5,-2.5)(2.5,2.5)
\psset{xunit=2.5,yunit=2.5}
\psline(-0.830,0.146)(-.866,-.500)
\psline(-0.830,0.146)(-.643,.766)
\psline(-.643,.766)(-0.502,0.421)
\psline(-.866,-.500)(-0.502,0.421)
\psline(.000,1.000)(-0.502,0.421)
\psline(-.866,-.500)(-0.074,0.421) 
\psline(.000,1.000)(-0.074,0.421) 
\psline(.643,.766)(-0.074,0.421) 
\psline(-.866,-.500)(0.254,0.146)
\psline(.643,.766)(0.254,0.146)
\psline(.985,.174)(0.254,0.146)
\psline(-.866,-.500)(-0.074,-0.422)
\psline(.985,.174)(-0.074,-0.422)
\psline(-.342,-.940)(-0.074,-0.422)
\psline(.985,.174)(0.328,-0.568)
\psline(-.342,-.940)(0.328,-0.568)
\psline(.342,-.940)(0.328,-0.568)
\psline(.985,.174)(0.730,-0.422)
\psline(.342,-.940)(0.730,-0.422)
\pscircle[fillstyle=solid,fillcolor=black](-.985,.174){.1}
\pscircle[fillstyle=solid,fillcolor=black](-.866,-.500){.2}
\pscircle[fillstyle=solid,fillcolor=black](-.643,.766){.2}
\pscircle[fillstyle=solid,fillcolor=black](.000,1.000){.2}
\pscircle[fillstyle=solid,fillcolor=black](.643,.766){.2}
\pscircle[fillstyle=solid,fillcolor=black](.985,.174){.2}
\pscircle[fillstyle=solid,fillcolor=black](-.342,-.940){.2}
\pscircle[fillstyle=solid,fillcolor=black](.342,-.940){.2}
\pscircle[fillstyle=solid,fillcolor=black](.866,-.500){.1}
\pscircle[fillstyle=solid,fillcolor=white](-0.830,0.146){.2}
\pscircle[fillstyle=solid,fillcolor=white](-0.502,0.421){.2}
\pscircle[fillstyle=solid,fillcolor=white](-0.074,0.421){.2} 
\pscircle[fillstyle=solid,fillcolor=white](0.254,0.146){.2}
\pscircle[fillstyle=solid,fillcolor=white](-0.074,-0.422){.2}
\pscircle[fillstyle=solid,fillcolor=white](0.328,-0.568){.2}
\pscircle[fillstyle=solid,fillcolor=white](0.730,-0.422){.2}
\endpspicture
$$
\begin{center}
Figure 5. A snake of order 6.
\end{center}

Although the operation of adding on a new 4-cycle, or ``box'',
can be done in three ways at each step 
in the iterative construction of a snake,
for purposes of counting perfect matchings
there are really only two choices at each stage:
for $k \geq 3$,
the edge that joins the $k-1$st box to the $k$th
can either (1) be disjoint from 
the edge that joins the $k-2$nd box to the $k-1$st
or (2) have an endpoint in common with it.
If we let $m_{k}$, $m_{k-1}$, and $m_{k-2}$
denote the number of perfect matchings of
the $k$th, $k-1$st, and $k-2$nd snakes
in the iterative process,
then (as we will now show)
in case 1, $m_{k} = m_{k-1} + m_{k-2}$,
while in case 2, $m_{k} = m_{k-1} + (m_{k-1} - m_{k-2})$;
that is, the three numbers 
are either in ``Fibonacci progression'' or in arithmetic progression.
Refer to Figure 6,
where $u$ and $v$ are the vertices of the $k-2$nd snake
that are not part of the $k-3$rd snake,
$w$ and $x$ are the vertices of the $k-1$st snake
that are not part of the $k-2$nd snake,
and $y$ and $z$ are the vertices of the $k$th snake
that are not part of the $k-1$st snake,
in the fashion shown.
In both cases, the difference $m_{k} - m_{k-1}$
counts those perfect matchings of the $k$th snake
that do not contain the edge $yz$
and that are therefore forced to contain
all the edges shown in bold in the figure.
In case 1, these perfect matchings
correspond to perfect matchings of the $k-2$nd snake:
simply delete edges $wy$ and $xz$ (and the vertices they contain).
In case 2, these perfect matchings
correspond to those perfect matchings of the $k-1$st snake
that do not arise from a perfect matching of the $k-2$nd snake
by adjoining the  edge $wx$:
simply delete edges $vy$ and $xz$ and adjoin edge $vx$.
$$
\pspicture(0.5,0.5)(7.5,4.5)
\rput(2,4){Case (1)}
\rput(6.5,4){Case (2)}
\psline(0.5,3)(0.8,3)
\psline(1.2,3)(1.8,3)
\psline[linewidth=3pt](2.2,3)(2.8,3)
\psline(0.5,2)(0.8,2)
\psline(1.2,2)(1.8,2)
\psline[linewidth=3pt](2.2,2)(2.8,2)
\psline(1,2.2)(1,2.8)
\psline(2,2.2)(2,2.8)
\psline(3,2.2)(3,2.8)
\psline(5.5,3)(5.8,3)
\psline[linewidth=3pt](6.2,3)(6.8,3)
\psline(5.5,2)(5.8,2)
\psline(6.2,2)(6.8,2)
\psline(6.2,1)(6.8,1)
\psline(6,2.2)(6,2.8)
\psline(7,2.2)(7,2.8)
\psline[linewidth=3pt](6,1.2)(6,1.8)
\psline[linewidth=3pt](7,1.2)(7,1.8)
\rput(1,3){$u$}
\rput(1,2){$v$}
\rput(2,3){$w$}
\rput(2,2){$x$}
\rput(3,3){$y$}
\rput(3,2){$z$}
\rput(6,3){$u$}
\rput(6,2){$v$}
\rput(7,3){$w$}
\rput(7,2){$x$}
\rput(6,1){$y$}
\rput(7,1){$z$}
\endpspicture
$$
\begin{center}
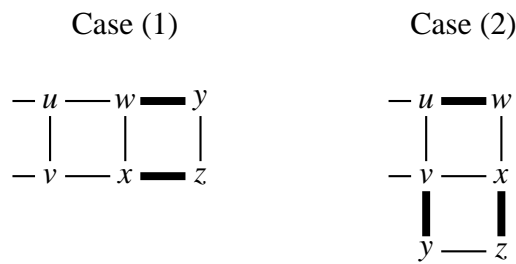
Figure 6. The snake-graph recurrence.
\end{center}

In terms of the triangulation picture,
a snake-graph of order $k$ corresponds to a chain of $k+1$ triangles,
in which the $i$th triangle (for $2 \leq i \leq k+1$)
consists of one edge $uv$ of the $i-1$st triangle
(not the edge joining the $i-1$st and $i-2$nd triangles)
together with a new vertex $x$
and two edges $ux$, $vx$.
Any two consecutive triangles in this chain share two vertices,
and any three consecutive triangles in this chain share one vertex.
If the last four consecutive triangles in the chain
have no vertex in common, then we are in case 1;
if they do have a vertex in common, we are in case 2.
We can make a {\it code\/} of length $n-2$
that contains this information.
Thus, the snake-graph shown in Figure 5
and the triangulation it arises from
can be described (from left to right) by the code
{\bf 2212}, indicating that as we move through the snake
from left to right,
we encounter case 2, case 2, case 1, and case 2.
Two extreme cases are the snake with code {\bf 11}...{\bf 1}
(the ``straight snake'')
and the snake with code {\bf 22}...{\bf 2}
(the ``fan snake'').

For purposes of enumeration of matchings,
every snake graph can be built as a chain of boxes
where each new box is either added at the right of the preceding box
or at the bottom of the preceding box.
This is because the two geometrically distinct subcases of case 2
are the same from the point of view of enumeration of matchings,
even though they are not isomorphic as graphs.
For instance, consider a triangulated $(n+3)$-gon
in which all the diagonals emanate from a single vertex.
Strictly speaking, the associated snake
(with code {\bf 22}...{\bf 2})
has all $n$ of its boxes sharing a single vertex.
However, we can replace this by a snake of
squares in the square grid,
where new squares are alternately added
at the right or at the bottom of the growing snake.
Both snake-graphs have exactly $n+1$ perfect matchings.
Similarly,
the snake-graph shown in Figure 5
and the snake-graph shown in Figure 7
both have the code {\bf 2212} and both have 13 perfect matchings.

$$
\pspicture(0.9,0.9)(5.1,4.1)
\psline(4,1)(5,1)
\psline(2,2)(5,2)
\psline(1,3)(5,3)
\psline(1,4)(3,4)
\psline(1,3)(1,4)
\psline(2,2)(2,4)
\psline(3,2)(3,4)
\psline(4,1)(4,3)
\psline(5,1)(5,3)
\pscircle[fillstyle=solid,fillcolor=black](1,4){.1}
\pscircle[fillstyle=solid,fillcolor=black](2,4){.1}
\pscircle[fillstyle=solid,fillcolor=black](3,4){.1}
\pscircle[fillstyle=solid,fillcolor=black](1,3){.1}
\pscircle[fillstyle=solid,fillcolor=black](2,3){.1}
\pscircle[fillstyle=solid,fillcolor=black](3,3){.1}
\pscircle[fillstyle=solid,fillcolor=black](4,3){.1}
\pscircle[fillstyle=solid,fillcolor=black](5,3){.1}
\pscircle[fillstyle=solid,fillcolor=black](2,2){.1}
\pscircle[fillstyle=solid,fillcolor=black](3,2){.1}
\pscircle[fillstyle=solid,fillcolor=black](4,2){.1}
\pscircle[fillstyle=solid,fillcolor=black](5,2){.1}
\pscircle[fillstyle=solid,fillcolor=black](4,1){.1}
\pscircle[fillstyle=solid,fillcolor=black](5,1){.1}
\endpspicture
$$
\begin{center}
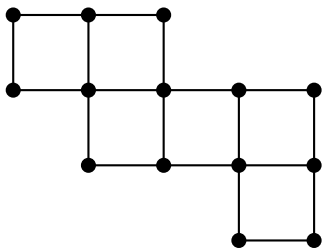
Figure 7. Another snake of order 6.
\end{center}

Graphs made of chains of hexagons have been considered before,
starting in the chemical literature
on account of their relevance to the study of benzenoid hydrocarbons;
an analogous theory applies there.
To add the $k$th hexagon to the chain,
we choose one of the three edges of the $k-1$st hexagon
that has no endpoints in common with the
edge that joins the $k-1$st and $k-2$nd hexagons in the chain.
If these two edges are diametrically opposite one another
in the hexagon that they both belong to,
we are in case 2, and the relation
$m_{k} - m_{k-1} = m_{k-1} - m_{k-2}$ applies;
otherwise, we are in case 1, and the relation
$m_{k} - m_{k-1} = m_{k-2}$ applies.

A good way to understand what is going on here
comes from consideration of products of the matrices
$A = \left( \begin{array}{cc} 0 & 1 \\ 1 & 1 \end{array} \right)$
and
$B = \left( \begin{array}{cc} 1 & 1 \\ 1 & 0 \end{array} \right)$.
A product of $n-1$ such matrices
corresponds to a snake with $n$ boxes,
where the presence of an $A$ (resp.\ $B$)
as the $i$th factor in the matrix product
(with $1 \leq i \leq n-1$)
corresponds to a horizontal (resp.\ vertical)
segment of the snake,
with the $i+1$st box in the snake
lying to the right of (resp.\ below)
the $i$th box in the snake.
For instance, the matrix product 
$ABAAB = \left( \begin{array}{cc} 2 & 1 \\ 7 & 3 \end{array} \right)$
corresponds to the snake-graph shown in Figure 7,
with code ${\bf 2212}$.
More generally, a product of $n-1$ matrices,
each of which is either $A$ or $B$,
corresponds to a snake of order $n$
whose code can be read off from the product by the following rule:
If the $i$th and $i+1$st matrices in the product are the same,
the $i$th element of the code is 1;
otherwise, the $i$th element of the code is 2.

The number of perfect matchings of a snake
is equal to the sum of the entries of the associated matrix
(so that in the specific example shown
the number of perfect matchings is $2+1+7+3$). 
More specifically: 
the matrix entry in the upper left counts the perfect matchings
in which both the upper-left vertex of the snake
and the lower-right vertex of the snake
are matched horizontally;
the matrix entry in the upper right counts the perfect matchings
in which the upper-left vertex is matched horizontally
and the lower-right vertex is matched vertically;
the matrix entry in the lower left counts the perfect matchings
in which the upper-left vertex is matched vertically
and the lower-right vertex is matched horizontally; and
the matrix entry in the lower right counts the perfect matchings
in which both the upper-left vertex 
and the lower-right vertex are matched vertically.
This interpretation of the entries of the product matrix
is easily verified by induction. 

$$
\pspicture(0.5,0.5)(6.5,2.5)
\psline(1,1)(2,2)
\psline(2,1)(3,2)
\psline(3,1)(4,2)
\psline(4,1)(5,2)
\psline(5,1)(6,2)
\psline(1,2)(2,1)
\psline(2,2)(3,1)
\psline(3,2)(4,1)
\psline(4,2)(5,1)
\psline(5,2)(6,1)
\psline(1,1)(2,1)
\psline(2,2)(3,2)
\psline(3,1)(5,1)
\psline(5,2)(6,2)
\pscircle[fillstyle=solid,fillcolor=black](1,1){.1}
\pscircle[fillstyle=solid,fillcolor=black](2,1){.1}
\pscircle[fillstyle=solid,fillcolor=black](3,1){.1}
\pscircle[fillstyle=solid,fillcolor=black](4,1){.1}
\pscircle[fillstyle=solid,fillcolor=black](5,1){.1}
\pscircle[fillstyle=solid,fillcolor=black](6,1){.1}
\pscircle[fillstyle=solid,fillcolor=black](1,2){.1}
\pscircle[fillstyle=solid,fillcolor=black](2,2){.1}
\pscircle[fillstyle=solid,fillcolor=black](3,2){.1}
\pscircle[fillstyle=solid,fillcolor=black](4,2){.1}
\pscircle[fillstyle=solid,fillcolor=black](5,2){.1}
\pscircle[fillstyle=solid,fillcolor=black](6,2){.1}
\endpspicture
$$
\begin{center}
Figure 8. A paths model for snakes.
\end{center}

A different combinatorial model	
that gives the same numbers as the perfect matchings model
arises from the fact that these numbers
can be expressed as the sum of the entries in a matrix
that is written as the product of matrices
whose entries are all 0's and 1's
(namely the matrices $A$ and $B$).
More specifically, we can create a graph
in which the number of paths from either of two source vertices
to either of two target vertices
is the same as the number of perfect matchings of a snake-graph.
Figure 8 shows what the paths-graph looks like 
for the snake associated with the matrix-product $ABAAB$.
Each factor of $A$ corresponds to a 4-vertex bipartite graph
containing all edges from the left to the right
except the edge connecting the top left to the top right,
and each factor of $B$ corresponds to a 4-vertex bipartite graph
containing all edges from the left to the right
except the edge connecting the bottom left to the bottom right.
Multiplication of matrices corresponds to adjunction of graphs,
and the definition of matrix multiplication
ensures that $i,j$th entry of the product
equals the number of left-right paths 
joining the $j$th of the two leftmost vertices
to the $i$th of the two rightmost vertices ($1 \leq i,j \leq 2$).
(For more on this combinatorial aspect of matrix multiplication,
see~\cite{Z1}.)
Note that changing all $A$'s into $B$'s and vice versa
simply flips the picture upside down.

We can improve on this model
by making a slight twist in our matrix-product,
working instead with the matrices
$L = \left( \begin{array}{cc} 1 & 1 \\ 0 & 1 \end{array} \right)$
and
$R = \left( \begin{array}{cc} 1 & 0 \\ 1 & 1 \end{array} \right)$.
To turn an $AB$-product into an essentially equivalent $LR$-product,
work from left to right,
with the proviso that two factors in the $LR$ product should be equal
if and only if two factors in the $AB$ product are {\it not\/}.
Thus, the product $ABAAB$
corresponds to either the product $LLLRR$ or the product $RRRLL$.
Either way, we get a product-matrix
whose four entries are, up to permutation, 
the same as the four entries of the $AB$ product,
with the virtue that the picture no longer involves crossings.
Figure 9, for instance, is the picture for $RRRLL$.

$$
\pspicture(0.5,0.5)(6.5,2.5)
\psline(1,2)(6,2)
\psline(1,1)(6,1)
\psline(1,2)(2,1)
\psline(2,2)(3,1)
\psline(3,2)(4,1)
\psline(4,1)(5,2)
\psline(5,1)(6,2)
\pscircle[fillstyle=solid,fillcolor=black](1,1){.1}
\pscircle[fillstyle=solid,fillcolor=black](2,1){.1}
\pscircle[fillstyle=solid,fillcolor=black](3,1){.1}
\pscircle[fillstyle=solid,fillcolor=black](4,1){.1}
\pscircle[fillstyle=solid,fillcolor=black](5,1){.1}
\pscircle[fillstyle=solid,fillcolor=black](6,1){.1}
\pscircle[fillstyle=solid,fillcolor=black](1,2){.1}
\pscircle[fillstyle=solid,fillcolor=black](2,2){.1}
\pscircle[fillstyle=solid,fillcolor=black](3,2){.1}
\pscircle[fillstyle=solid,fillcolor=black](4,2){.1}
\pscircle[fillstyle=solid,fillcolor=black](5,2){.1}
\pscircle[fillstyle=solid,fillcolor=black](6,2){.1}
\endpspicture
$$
\begin{center}
Figure 9. A planar paths model for snakes.
\end{center}

A variant of this picture is shown in Figure 10.
This is just like the Figure 9,
except that we have added a vertex at the upper left
that connects to the two previously leftmost vertices,
and we have added a vertex at the lower right
that connects to the two previously rightmost vertices,
so that, where before we counted paths
from either of the two leftmost vertices
to either of the two rightmost vertices
(obtaining four numbers that get added together),
we now count paths
from the unique leftmost vertex
to the unique rightmost vertex.
We have marked each vertex $v$ by a number
that indicates the number of paths
from the leftmost vertex to $v$.
The leftmost vertex gets marked with a 1,
and every other vertex gets marked with the sum
of the numbers marking its (one or two) leftward neighbors.
In terms of the triangulation,
this means that we put 1's
at the vertices of the initial triangle in the snake,
and we proceed marking vertices along the snake all the way to its tail,
where each new vertex is marked with the sum
of the markings of the other two vertices
of the triangle being added to the snake.

$$
\pspicture(-0.5,0.5)(7.5,2.5)
\psline(0,2)(6,2)
\psline(1,1)(7,1)
\psline(0,2)(1,1)
\psline(1,2)(2,1)
\psline(2,2)(3,1)
\psline(3,2)(4,1)
\psline(4,1)(5,2)
\psline(5,1)(6,2)
\psline(6,2)(7,1)
\pscircle[fillstyle=solid,fillcolor=black](1,1){.1}
\pscircle[fillstyle=solid,fillcolor=black](2,1){.1}
\pscircle[fillstyle=solid,fillcolor=black](3,1){.1}
\pscircle[fillstyle=solid,fillcolor=black](4,1){.1}
\pscircle[fillstyle=solid,fillcolor=black](5,1){.1}
\pscircle[fillstyle=solid,fillcolor=black](6,1){.1}
\pscircle[fillstyle=solid,fillcolor=black](7,1){.1}
\pscircle[fillstyle=solid,fillcolor=black](0,2){.1}
\pscircle[fillstyle=solid,fillcolor=black](1,2){.1}
\pscircle[fillstyle=solid,fillcolor=black](2,2){.1}
\pscircle[fillstyle=solid,fillcolor=black](3,2){.1}
\pscircle[fillstyle=solid,fillcolor=black](4,2){.1}
\pscircle[fillstyle=solid,fillcolor=black](5,2){.1}
\pscircle[fillstyle=solid,fillcolor=black](6,2){.1}
\rput(1,0.6){1}
\rput(2,0.6){2}
\rput(3,0.6){3}
\rput(4,0.6){4}
\rput(5,0.6){4}
\rput(6,0.6){4}
\rput(7,0.6){13}
\rput(0,2.4){1}
\rput(1,2.4){1}
\rput(2,2.4){1}
\rput(3,2.4){1}
\rput(4,2.4){1}
\rput(5,2.4){5}
\rput(6,2.4){9}
\endpspicture
$$
\begin{center}
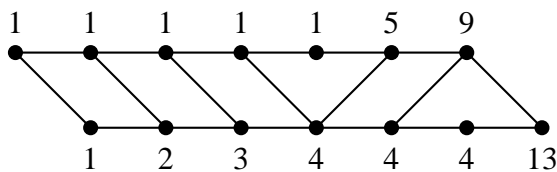
Figure 10. Another planar paths model for snakes.
\end{center}

The marking scheme of Figure 10 is in fact nothing more than
a slight variation on
Conway and Coxeter's method of computing entries in frieze patterns
by successively marking the vertices in a triangulation.
Conway and Coxeter mark a single vertex with a 0, all its neighbors with 1's,
and proceed from there;
since we are pruning away all the side-branches
of the graph $G_{i,j}$ until all that remains is a snake,
we are effectively limiting ourselves to the case
where the vertex to be marked with a 0 has only two neighbors.
In this case, the only difference
between our marking scheme and Conway and Coxeter's
is that they mark the initial vertex with a 0
where we mark it with a 1.
Figure 11 shows what Conway and Coxeter's scheme looks like for the snake
whose different representations are shown in Figures 5 through 10.
The 9-gon being triangulated is not convex, 
but that is not a problem since we are
dealing with triangulations purely combinatorially.

$$
\psset{xunit=0.7,yunit=0.7}
\pspicture(0,1.366)(6,7.294)
\psline(3,1.732)(5,1.732)
\psline(2,3.464)(4,3.464)
\psline(1,5.196)(5,5.196)
\psline(2,6.928)(4,6.928)
\psline(3,1.732)(5,5.196)
\psline(2,3.464)(4,6.928)
\psline(1,5.196)(2,6.928)
\psline(2,3.464)(3,1.732)
\psline(2,6.928)(5,1.732)
\psline(4,6.928)(5,5.196)
\rput(0.8,5.196){0}
\rput(1.8,6.928){1}
\rput(2.7,4.996){1}
\rput(4.2,6.928){2}
\rput(5.2,5.196){3}
\rput(4.2,3.464){4}
\rput(1.8,3.464){5}
\rput(2.8,1.732){9}
\rput(5.2,1.732){13}
\endpspicture
$$
\begin{center}
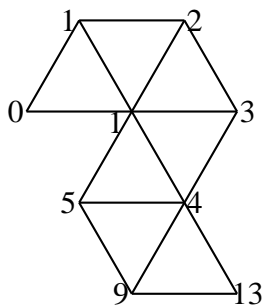
Figure 11. Conway and Coxeter's marking scheme.
\end{center}

At this point it should be mentioned that there is a link
between the directed path model of Figure 10,
the hexagon snake model, 
and the square snake model,
by way of a multigraph matching model
that is in turn related
to the strip-tiling model of Benjamin and Quinn~\cite{BQ}.
We start by making use of a correspondence
that has been rediscovered a number of times,
starting as far back as 1952 \cite{GD} \cite{S} \cite{Kup}:

Proposition: Let $D$ be a directed acyclic graph with vertex set $V$,
where $m$ vertices $s_1,\dots,s_m$ have been designated as sources
and $m$ other vertices $t_1,\dots,t_m$ have been designated as targets.
(Since $D$ is acyclic, self-loops are forbidden,
but $D$ is permitted to have multiple edges.)
Create an undirected graph $G$ with vertex set $V \times \{1,2\}$
and two kinds of edges:
for each vertex $v$ of $D$, 
$D'$ contains an edge joining $(v,1)$ and $(v,2)$,
and for each directed edge $e: v \rightarrow w$ of $D$,
$D'$ contains an edge joining $(v,2)$ and $(w,1)$.
(If $D$ has more than one directed edge from $v$ to $w$,
$G$ has just as many edges joining $(v,2)$ and $(w,1)$.)
Let $H$ be the induced subgraph of $G$
obtained by deleting all the vertices $(s_i,1)$ and $(t_i,2)$
($1 \leq i \leq m$)
and all incident edges.
Then the perfect matchings of $H$
are equinumerous with the ways to join
the $m$ sources to the $m$ targets by $m$ edge-disjoint paths in $D$
(which need not connect $s_1$ to $t_1$ etc.).
Specifically, given such a collection of $m$ paths,
take each arc $e: v \rightarrow w$ that belongs to one of the paths
and replace it by the corresponding edge joining $(v,2)$ and $(w,1)$ in $H$,
and replace each vertex $v$ in $D$ that does not lie on any of the paths
by the edge joining $(v,1)$ and $(v,2)$ in $H$.  
It is easy to check that this yields a perfect matching of $H$,
and it is also easy to show that every perfect matching
arises in a unique fashion in this way.

If we carry out this operation with the directed graph $D$
shown in Figure 10
(where all edges are oriented from left to right
with the leftmost vertex the sole source
and the rightmost vertex the sole target),
we obtain a graph $G$ composed of 6-cycles (hexagons),
as shown in Figure 12.
$$
\psset{xunit=1.0,yunit=1.0}
\pspicture(-0.5,0.5)(7.5,2.5)
\psline(0,2)(6.2,2)
\psline(0.8,1)(6.8,1)
\psline(0,2)(0.8,1)
\psline(1.2,2)(1.8,1)
\psline(2.2,2)(2.8,1)
\psline(3.2,2)(3.8,1)
\psline(4.2,1)(4.8,2)
\psline(5.2,1)(5.8,2)
\psline(6.2,2)(6.8,1)
\pscircle[fillstyle=solid,fillcolor=black](0.8,1){.1}
\pscircle[fillstyle=solid,fillcolor=black](1.2,1){.1}
\pscircle[fillstyle=solid,fillcolor=black](1.8,1){.1}
\pscircle[fillstyle=solid,fillcolor=black](2.2,1){.1}
\pscircle[fillstyle=solid,fillcolor=black](2.8,1){.1}
\pscircle[fillstyle=solid,fillcolor=black](3.2,1){.1}
\pscircle[fillstyle=solid,fillcolor=black](3.8,1){.1}
\pscircle[fillstyle=solid,fillcolor=black](4.2,1){.1}
\pscircle[fillstyle=solid,fillcolor=black](4.8,1){.1}
\pscircle[fillstyle=solid,fillcolor=black](5.2,1){.1}
\pscircle[fillstyle=solid,fillcolor=black](5.8,1){.1}
\pscircle[fillstyle=solid,fillcolor=black](6.2,1){.1}
\pscircle[fillstyle=solid,fillcolor=black](6.8,1){.1}
\pscircle[fillstyle=solid,fillcolor=black](0,2){.1}
\pscircle[fillstyle=solid,fillcolor=black](0.8,2){.1}
\pscircle[fillstyle=solid,fillcolor=black](1.2,2){.1}
\pscircle[fillstyle=solid,fillcolor=black](1.8,2){.1}
\pscircle[fillstyle=solid,fillcolor=black](2.2,2){.1}
\pscircle[fillstyle=solid,fillcolor=black](2.8,2){.1}
\pscircle[fillstyle=solid,fillcolor=black](3.2,2){.1}
\pscircle[fillstyle=solid,fillcolor=black](3.8,2){.1}
\pscircle[fillstyle=solid,fillcolor=black](4.2,2){.1}
\pscircle[fillstyle=solid,fillcolor=black](4.8,2){.1}
\pscircle[fillstyle=solid,fillcolor=black](5.2,2){.1}
\pscircle[fillstyle=solid,fillcolor=black](5.8,2){.1}
\pscircle[fillstyle=solid,fillcolor=black](6.2,2){.1}
\endpspicture
$$
\begin{center}
Figure 12. A hexagon snake.
\end{center}

\noindent
Figure 13 shows how a particular path in $D$ 
corresponds to a particular perfect matching in $G$
(as described by the above proof).

$$
\pspicture(-0.5,0.5)(7.5,5.5)
\psline(0,5)(6,5)
\psline(1,4)(7,4)
\psline(0,5)(1,4)
\psline(1,5)(2,4)
\psline(2,5)(3,4)
\psline(3,5)(4,4)
\psline(4,4)(5,5)
\psline(5,4)(6,5)
\psline(6,5)(7,4)
\psline[linewidth=3pt](0,5)(2,5)
\psline[linewidth=3pt](2,5)(3,4)
\psline[linewidth=3pt](3,4)(5,4)
\psline[linewidth=3pt](5,4)(6,5)
\psline[linewidth=3pt](6,5)(7,4)
\pscircle[fillstyle=solid,fillcolor=black](1,4){.1}
\pscircle[fillstyle=solid,fillcolor=black](2,4){.1}
\pscircle[fillstyle=solid,fillcolor=black](3,4){.1}
\pscircle[fillstyle=solid,fillcolor=black](4,4){.1}
\pscircle[fillstyle=solid,fillcolor=black](5,4){.1}
\pscircle[fillstyle=solid,fillcolor=black](6,4){.1}
\pscircle[fillstyle=solid,fillcolor=black](7,4){.1}
\pscircle[fillstyle=solid,fillcolor=black](0,5){.1}
\pscircle[fillstyle=solid,fillcolor=black](1,5){.1}
\pscircle[fillstyle=solid,fillcolor=black](2,5){.1}
\pscircle[fillstyle=solid,fillcolor=black](3,5){.1}
\pscircle[fillstyle=solid,fillcolor=black](4,5){.1}
\pscircle[fillstyle=solid,fillcolor=black](5,5){.1}
\pscircle[fillstyle=solid,fillcolor=black](6,5){.1}
\rput(3.5,3){\Large{$\downarrow$}}
\psline(0.2,2)(6.2,2)
\psline(0.8,1)(6.8,1)
\psline(0.2,2)(0.8,1)
\psline(1.2,2)(1.8,1)
\psline(2.2,2)(2.8,1)
\psline(3.2,2)(3.8,1)
\psline(4.2,1)(4.8,2)
\psline(5.2,1)(5.8,2)
\psline(6.2,2)(6.8,1)
\psline[linewidth=3pt](0.8,1)(1.2,1)
\psline[linewidth=3pt](1.8,1)(2.2,1)
\psline[linewidth=3pt](3.2,1)(3.8,1)
\psline[linewidth=3pt](4.2,1)(4.8,1)
\psline[linewidth=3pt](5.8,1)(6.2,1)
\psline[linewidth=3pt](0.2,2)(0.8,2)
\psline[linewidth=3pt](1.2,2)(1.8,2)
\psline[linewidth=3pt](2.8,2)(3.2,2)
\psline[linewidth=3pt](3.8,2)(4.2,2)
\psline[linewidth=3pt](4.8,2)(5.2,2)
\psline[linewidth=3pt](2.2,2)(2.8,1)
\psline[linewidth=3pt](5.2,1)(5.8,2)
\psline[linewidth=3pt](6.2,2)(6.8,1)
\pscircle[fillstyle=solid,fillcolor=black](0.8,1){.1}
\pscircle[fillstyle=solid,fillcolor=black](1.2,1){.1}
\pscircle[fillstyle=solid,fillcolor=black](1.8,1){.1}
\pscircle[fillstyle=solid,fillcolor=black](2.2,1){.1}
\pscircle[fillstyle=solid,fillcolor=black](2.8,1){.1}
\pscircle[fillstyle=solid,fillcolor=black](3.2,1){.1}
\pscircle[fillstyle=solid,fillcolor=black](3.8,1){.1}
\pscircle[fillstyle=solid,fillcolor=black](4.2,1){.1}
\pscircle[fillstyle=solid,fillcolor=black](4.8,1){.1}
\pscircle[fillstyle=solid,fillcolor=black](5.2,1){.1}
\pscircle[fillstyle=solid,fillcolor=black](5.8,1){.1}
\pscircle[fillstyle=solid,fillcolor=black](6.2,1){.1}
\pscircle[fillstyle=solid,fillcolor=black](6.8,1){.1}
\pscircle[fillstyle=solid,fillcolor=black](0.2,2){.1}
\pscircle[fillstyle=solid,fillcolor=black](0.8,2){.1}
\pscircle[fillstyle=solid,fillcolor=black](1.2,2){.1}
\pscircle[fillstyle=solid,fillcolor=black](1.8,2){.1}
\pscircle[fillstyle=solid,fillcolor=black](2.2,2){.1}
\pscircle[fillstyle=solid,fillcolor=black](2.8,2){.1}
\pscircle[fillstyle=solid,fillcolor=black](3.2,2){.1}
\pscircle[fillstyle=solid,fillcolor=black](3.8,2){.1}
\pscircle[fillstyle=solid,fillcolor=black](4.2,2){.1}
\pscircle[fillstyle=solid,fillcolor=black](4.8,2){.1}
\pscircle[fillstyle=solid,fillcolor=black](5.2,2){.1}
\pscircle[fillstyle=solid,fillcolor=black](5.8,2){.1}
\pscircle[fillstyle=solid,fillcolor=black](6.2,2){.1}
\endpspicture
$$
\begin{center}
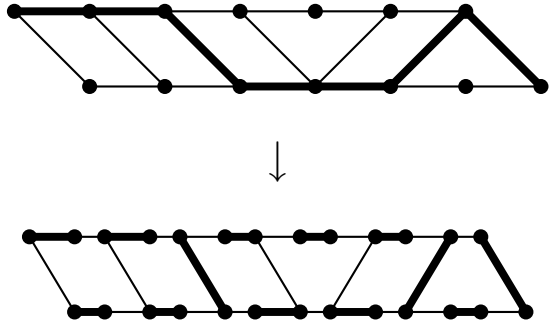
Figure 13. From path systems to perfect matchings.
\end{center}

The hexagon snake of Figure 12
can be drawn as a snake of regular hexagons,
as shown in Figure 14.
If we turn the figure on its side,
so that the leftmost square is at the top,
we can see how the turns of the snake
correspond to the symbols in its $L,R$-string $RRRLL$.

$$
\psset{xunit=0.25,yunit=0.433}
\pspicture(-1.5,1.5)(21,6.5)
\psline(1,3)(2,2)
\psline(1,3)(2,4)
\psline(2,2)(4,2)
\psline(2,4)(4,4)
\psline(4,2)(5,3)
\psline(4,4)(5,3)
\psline(4,4)(5,5)
\psline(5,3)(7,3)
\psline(5,5)(7,5)
\psline(7,3)(8,4)
\psline(7,5)(8,4)
\psline(7,5)(8,6)
\psline(8,4)(10,4)
\psline(8,6)(10,6)
\psline(10,4)(11,5)
\psline(10,6)(11,5)
\psline(10,6)(11,7)
\psline(11,5)(13,5)
\psline(11,7)(13,7)
\psline(13,5)(14,4)
\psline(13,5)(14,6)
\psline(13,7)(14,6)
\psline(14,4)(16,4)
\psline(14,6)(16,6)
\psline(16,4)(17,3)
\psline(16,4)(17,5)
\psline(16,6)(17,5)
\psline(17,3)(19,3)
\psline(17,5)(19,5)
\psline(19,3)(20,4)
\psline(19,5)(20,4)
\pscircle[fillstyle=solid,fillcolor=black](1,3){.1}
\pscircle[fillstyle=solid,fillcolor=black](2,2){.1}
\pscircle[fillstyle=solid,fillcolor=black](2,4){.1}
\pscircle[fillstyle=solid,fillcolor=black](4,2){.1}
\pscircle[fillstyle=solid,fillcolor=black](4,4){.1}
\pscircle[fillstyle=solid,fillcolor=black](5,3){.1}
\pscircle[fillstyle=solid,fillcolor=black](5,5){.1}
\pscircle[fillstyle=solid,fillcolor=black](7,3){.1}
\pscircle[fillstyle=solid,fillcolor=black](7,5){.1}
\pscircle[fillstyle=solid,fillcolor=black](8,4){.1}
\pscircle[fillstyle=solid,fillcolor=black](8,6){.1}
\pscircle[fillstyle=solid,fillcolor=black](10,4){.1}
\pscircle[fillstyle=solid,fillcolor=black](10,6){.1}
\pscircle[fillstyle=solid,fillcolor=black](11,5){.1}
\pscircle[fillstyle=solid,fillcolor=black](11,7){.1}
\pscircle[fillstyle=solid,fillcolor=black](13,5){.1}
\pscircle[fillstyle=solid,fillcolor=black](13,7){.1}
\pscircle[fillstyle=solid,fillcolor=black](14,4){.1}
\pscircle[fillstyle=solid,fillcolor=black](14,6){.1}
\pscircle[fillstyle=solid,fillcolor=black](16,4){.1}
\pscircle[fillstyle=solid,fillcolor=black](16,6){.1}
\pscircle[fillstyle=solid,fillcolor=black](17,3){.1}
\pscircle[fillstyle=solid,fillcolor=black](17,5){.1}
\pscircle[fillstyle=solid,fillcolor=black](19,3){.1}
\pscircle[fillstyle=solid,fillcolor=black](19,5){.1}
\pscircle[fillstyle=solid,fillcolor=black](20,4){.1}
\endpspicture
$$
\begin{center}
Figure 14. A hexagon snake with regular hexagons.
\end{center}

We also make use of an even simpler proposition
that is part of the folklore of perfect matchings:
Suppose $v$ is a vertex of degree 2 in a graph $G$,
with neighbors $w_1$ and $w_2$.
Let $G'$ be the graph obtained from $G$
by deleting $v$ and its two edges
and identifying vertices $w_1$ and $w_2$,
so that the new amalgamated vertex (call it $w$)
inherits the neighbors of $w_1$ and $w_2$.
(Specifically, if $x$ is some vertex
that in $G$ is connected to $w_1$ by $m_1$ edges
and connected to $w_2$ by $m_2$ edges,
then in $G'$, $x$ is connected to $w$ by $m_1+m_2$ edges.)
Then the perfect matchings of $G$
are equinumerous with the perfect matchings of $G'$.
Specifically, given a perfect matching of $G$
in which $v$ is connected to one of $w_1,w_2$
and the other is connected a vertex $x$,
construct a perfect matching of $G'$
in which $w$ is connected to $x$
and all other edges are unaffected
(in the case where there are multiple edges from $w$ to $x$,
one uses the edge that is associated with
the specific edge of the matching of $G$ that contains $x$).

Using this path-contraction operation,
one can show that enumeration of perfect matchings
of an arbitrary snake formed from $n$ cycles of even order
(i.e., any combination of 4-cycles, 6-cycles, etc., arranged in a chain
consisting of $n$ cycles)
reduces to enumeration of perfect matchings
of an ordinary straight snakes (made of 4-cycles)
in which the edges shared by one 4-cycle and the next
are allowed to have multiplicity,
with multiplicities adding up to $n+1$.
These multiplicities can be easily read off from
the $L,R$-string associated with the snake:
simply duplicate the first and last symbols of the string,
and then write down the run-lengths.
For instance, the $L,R$-string $RRRLL$
becomes $RRRRLLL$ when the first and last symbols are duplicated,
which gives the sequence of run-lengths $4,3$,
so that the graphs shown in Figures 7 and 13,
when contracted,
both become the multigraph shown in Figure 15,
where the 4 represents 4 parallel edges
and the 3 represents 3 parallel edges.

$$
\psset{xunit=0.7,yunit=0.7}
\pspicture(0.5,0.5)(3.5,3.5)
\psline(1,1)(3,1)
\psline(1,3)(3,3)
\psline(1,1)(1,3)
\psline(3,1)(3,3)
\rput(0.8,2){4}
\rput(3.2,2){3}
\pscircle[fillstyle=solid,fillcolor=black](1,1){.1}
\pscircle[fillstyle=solid,fillcolor=black](1,3){.1}
\pscircle[fillstyle=solid,fillcolor=black](3,1){.1}
\pscircle[fillstyle=solid,fillcolor=black](3,3){.1}
\endpspicture
$$
\begin{center}
Figure 15. A (short) straight snake with multiplicities.
\end{center}

The perfect matchings of such a weighted graph
can in turn be associated with strip-tilings
of the sort considered by Benjamin and Quinn~\cite{BQ}.
Specifically, suppose we have a straight snake of order $n$
whose $n+1$ vertical edges have multiplicities $r_0,r_1,\dots,r_N$.
Then we associate this with a 1-by-$(n+1)$ rectangular strip
that is to be covered by stackable 1-by-1 square tiles
and non-stackable 1-by-2 rectangular tiles (``dominos''),
where each square in the strip must be covered by at least one tile,
and where square tiles may be stacked to height $r_i$
at the $i$th square of the strip.
E.g., for the graph shown in Figure 15,
the associated strip-tiling problem
would involve a strip consisting of two squares,
which can either be tiled by a single domino
or by two non-empty stacks of squares
(up to 4 squares in the left stack
and up to 3 squares in the right stack).

By this point in the article, many readers will have recognized
that our combinatorial model is not too far removed
from the theory of continued fractions.
Work of Benjamin and Quinn, in the context of the strip-tiling model,
shows how combinatorial models 
can illuminate facts about continued fractions
(especially those like~\cite{BQS} and~\cite{BZ} that involve 
reversing the order of the convergents:
this operation seems unnatural from the point of view
of the definition of continued fractions,
inasmuch as it switches the high-order and low-order parts
of the continued fraction representation,
but the operation is extremely natural
for tilings of a strip).


There is a different way to relate frieze patterns to snake-graphs,
where we count paths in the snake-graphs themselves.
For instance, the number 13,
whose various enumerative interpretations 
we have followed throughout this section,
also occurs as the number of paths 
from the leftmost vertex to the rightmost vertex 
in the hexagon snake shown in Figure 14.
To see why, note that for purposes of enumerating such paths,
we can shrink each horizontal edge in Figure 14 to a point
(identifying the two endpoints),
obtaining a square snake (see Figure 16)
that is combinatorially the same as
the square snake shown in Figure 10.
It should be stressed that this square snake is not the square snake
we started with (shown in Figure 7).
It is ``dual'' to the original square snake,
making a bend where the original snake goes straight
and going straight where the original snake makes a bend.
(Equivalently, the code of the first snake has a {\bf 1}
where the code of the second snake has a {\bf 2},
and vice versa.)
Enumerating perfect matchings of each snake
is equivalent to counting paths its dual (from head to tail).
For instance, the snake in Figure 7 has 13 perfect matchings
and 19 paths from head to tail,
while the snake in Figure 16 has 19 perfect matchings
and 13 paths from head to tail.

$$
\psset{xunit=0.7,yunit=0.7}
\pspicture(1,1)(8,6)
\psline(1,2)(2,1)
\psline(1,2)(2,3)
\psline(2,1)(3,2)
\psline(2,3)(3,2)
\psline(2,3)(3,4)
\psline(3,2)(4,3)
\psline(3,4)(4,3)
\psline(3,4)(4,5)
\psline(4,3)(5,4)
\psline(4,5)(5,4)
\psline(4,5)(5,6)
\psline(5,4)(6,3)
\psline(5,4)(6,5)
\psline(5,6)(6,5)
\psline(6,3)(7,2)
\psline(6,3)(7,4)
\psline(6,5)(7,4)
\psline(7,2)(8,3)
\psline(7,4)(8,3)
\endpspicture
$$
\begin{center}
Figure 16. A dual snake.
\end{center}

$$
\psset{xunit=0.7,yunit=0.7}
\pspicture(0,0)(9,9)
\rput(4,1){$1_d$}
\rput(3,2){$1_c$}
\rput(5,2){4}
\rput(2,3){1}
\rput(4,3){3}
\rput(6,3){7}
\rput(1,4){1}
\rput(3,4){2}
\rput(5,4){5}
\rput(7,4){10}
\rput(2,5){1}
\rput(4,5){3}
\rput(6,5){7}
\rput(8,5){13}
\rput(3,6){1}
\rput(5,6){4}
\rput(7,6){9}
\rput(4,7){$1_b$}
\rput(6,7){5}
\rput(5,8){$1_a$}
\psline(3.8,1.2)(3.2,1.8)
\psline(2.8,2.2)(2.2,2.8)
\psline(1.8,3.2)(1.2,3.8)
\psline(4.8,2.2)(4.2,2.8)
\psline(3.8,3.2)(3.2,3.8)
\psline(2.8,4.2)(2.2,4.8)
\psline(3.8,5.2)(3.2,5.8)
\psline(4.8,6.2)(4.2,6.8)
\psline(5.8,7.2)(5.2,7.8)
\psline(4.2,1.2)(4.8,1.8)
\psline(3.2,2.2)(3.8,2.8)
\psline(2.2,3.2)(2.8,3.8)
\psline(3.2,4.2)(3.8,4.8)
\psline(4.2,5.2)(4.8,5.8)
\psline(5.2,6.2)(5.8,6.8)
\psline(1.2,4.2)(1.8,4.8)
\psline(2.2,5.2)(2.8,5.8)
\psline(3.2,6.2)(3.8,6.8)
\psline(4.2,7.2)(4.8,7.8)
\endpspicture
$$
\begin{center}
Figure 17. A dual snake in a frieze pattern.
\end{center}

There is a nice way to see a dual square snake of order $n$
as residing within a frieze pattern of order $n+3$:
rotate the snake by 90 degrees,
so that its first cell is at the top
and its last cell is at the bottom,
and put its top vertex (call it $u$) in the initial row of 1's
of an initially blank frieze pattern with $n+2$ rows,
so that its bottom vertex (call it $v$) lands in the final row of 1's,
and the $R$'s and $L$'s indicate
whether each successive box in the snake
lies to the right or left (respectively)
of the previous box.
If we put 1's along the left border of the snake,
we get a zig-zag of the kind discussed earlier,
so we obtain a frieze-pattern.
Moreover, within the part of the frieze-pattern
that is bounded by the line of slope $-1$ through $u$,
the line of slope $+1$ through $v$,
and the snake itself,
each entry admits an enumerative interpretation
relating to paths in the snake graph.
Specifically, given any location $w$ in the table
in the aforementioned region,
let $u'$ be the leftmost place 
where the line through $w$ of slope $-1$ meets the snake,
and let $v'$ be the leftmost place 
where the line through $w$ of slope $+1$ meets the snake;
then the entry at $w$ is equal to
the number of paths in the snake
from $u'$ to $v'$.
For instance, in Figure 17
which shows what happens for the $LLLRR$ snake
(with subscripts attached to some of the 1's for purposes of labeling), 
we find that 
9 is the number of paths from $1_a$ to $1_c$,
7 is the number of paths from $1_b$ to $1_c$,
13 is the number of paths from $1_a$ to $1_d$, and
10 is the number of paths from $1_b$ to $1_d$.
(The reader may find it instructive to compare this picture 
with the corresponding picture for the $RRRLL$ snake;
the arithmetic calculations are different,
but the number 13 still emerges as the rightmost entry.)

To see why this connection between dual snakes and frieze patterns holds,
we can use Lindstr\"om's lemma~\cite{L}
(rediscovered later by John and Sachs~\cite{JS}
and by Gessel and Viennot~\cite{GV}).
The $m=2$ case of this lemma states that 
if we have a directed graph 
with sources $s_1,s_2$ and targets $t_1,t_2$,
and there is no way to create a pair of vertex-disjoint paths
that join $s_1$ to $t_2$ and $s_2$ to $t_1$ respectively,
then the number of ways to create a pair of vertex-disjoint paths
that join $s_1$ to $t_1$ and $s_2$ to $t_2$ respectively
is equal to the 2-by-2 determinant $p_{11} p_{22} - p_{12} p_{22}$,
where $p_{ij}$ denotes the number of paths
from $s_i$ to $t_j$.
In our example, putting 
$s_1 = 1_a$, $s_2 = 1_b$, $t_1 = 1_d$, and $t_2 = 1_c$ [sic],
we see that there is no way to create
a path from $1_a$ to $1_c$
and a path from $1_b$ to $1_d$
that do not cross,
so the hypotheses are satisfied.
Furthermore, there is exactly 1 way
to create a path from $1_a$ to $1_d$
and a path from $1_b$ to $1_c$
that do not cross,
so we may conclude that
$p_{11} p_{22} - p_{12} p_{22}$ equals 1,
which (given that the $p_{ij}$'s are entries in our table)
is exactly the frieze relation.

We mention two other combinatorial models of frieze patterns,
for the sake of completeness:
Gregory Price's paths model~\cite{CP}
and the model of Broline, Crowe and Isaacs~\cite{BCI}.
The former (which has been significantly generalized
by Schiffler and Thomas; see~\cite{ST})
is related to the perfect matching model
by the bijection of Carroll and Price,
and the latter is closely related
to the Conway-Coxeter marking scheme.
      
\section{A tropical analogue} \label{sec-tropical}

Since the sideways frieze relation
involves only subtraction-free expressions in the cluster variables,
our whole picture admits a tropical analogue
(for background on tropical mathematics, see~\cite{SS})
in which multiplication is replaced by addition,
division by subtraction, addition by max, and 1 by 0.
(One could use min instead of max, but max will be more useful for us.)
In this new picture, the Ptolemy relation
$$d_{i,j} \, d_{k,l} + d_{j,k} \, d_{i,l} = d_{i,k} \, d_{j,l}$$
becomes the ultrametric relation
$$\max(d_{i,j}+d_{k,l}, d_{j,k}+d_{i,l}) = d_{i,k}+d_{j,l}.$$
Metrics satisfying this relation
arise from finite collections of non-intersecting arcs
that join points on the sides of the $n$-gon in pairs,
where the endpoints of such an arc are not permitted
to be vertices of the $n$-gon.
We will call such a collection of arcs 
an {\it integral lamination\/}.
Figure 18 shows an integral lamination
of a hexagon.

$$
\pspicture(0,-.75)(9,5.796)
\psset{xunit=1.5,yunit=1.5}
\pspolygon(2,0)(4,0)(5,1.732)(4,3.464)(2,3.464)(1,1.732)
\psarc(1,1.732){.8}{300}{60}
\psarc(2,3.464){.8}{240}{360}
\psarc(4.5,.866){.7}{60}{240}
\pscurve(1.5,2.598)(2.4,1.299)(2.8,0)
\pscurve(4.5,2.598)(3.6,1.299)(3.2,0)
\endpspicture
$$
\begin{center}
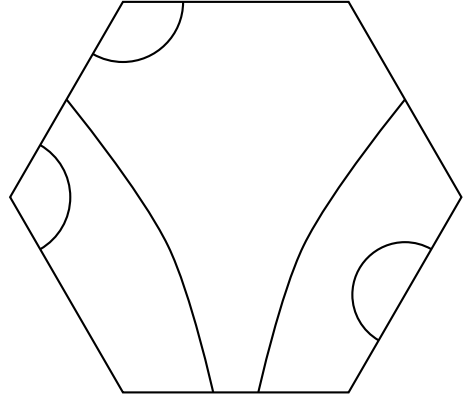
Figure 18.  An integral lamination.
\end{center}

\noindent
For any pair of vertices $i,j$,
we define $d_{i,j}$ as the smallest possible number of intersections
between a path in the $n$-gon from $i$ to $j$
and the arcs in the integral lamination
(we choose the path so as to avoid crossing any arc 
in the integral lamination more than once).
Then these quantities $d_{i,j}$ satisfy the ultrametric relation,
and thus can be arranged to form a tropical frieze pattern 
satisfying the relation
$$
\begin{array}{ccccccccl}
 X  &     &     &     &     &     &  Y  &   & \\
    & \jd &     &  A  &     & \ji &     &   & \\[1.2ex]
    &     &  B  &     &  C  &     &     & : & \ \ B + C = \max(A+D,X+Y) \\
    & \ji &     &  D  &     & \jd &     &   & \\
 Y  &     &     &     &     &     &  X  &   & \\
\end{array}
$$
For instance, the integral lamination of Figure 18
gives rise to the tropical frieze pattern 
$$
\begin{array}{ccccccccccccccc}
...&3& &1& &1& &0& &2& &1& &3&... \\[1.1ex]
...& &2& &2& &1& &2& &3& &2& &... \\[1.1ex]
...&1& &3& &2& &1& &3& &2& &1&... \\[1.1ex]
...& &2& &3& &2& &2& &2& &1& &... \\[1.1ex]
...&0& &2& &1& &3& &1& &1& &0&...
\end{array}
$$
As in the non-tropical case, we can find all the quantities $d_{i,j}$
once we know the values for all $(i,j)$ associated with
with the sides and diagonals 
belonging to some triangulation of the $n$-gon.

For an alternative picture,
one can divide the laminated $n$-gon into a finite number of sub-regions,
each of which is bounded by pieces of the boundary of the $n$-gon
and/or arcs of the integral lamination;
the vertices of the $n$-gon correspond to $n$ special sub-regions
(some of which may coincide with one another, 
if there is no arc in the integral lamination
separating the associated vertices of the $n$-gon).
Then the dual of this dissection of the $n$-gon
is a tree with $n$ specified leaf vertices (some of which may coincide),
and $d_{i,j}$ is the graph-theoretic distance
between leaf $i$ and leaf $j$ (which could be zero).
We see that if we know $2n-3$ of these leaf-to-leaf distances,
and the $2n-3$ pairs of leaves correspond
to the sides and diagonals of a triangulated $n$-gon,
then all of the other leaf-to-leaf distances
can be expressed as piecewise-linear functions
(involving just plus, minus, and max)
of the $2n-3$ specified distances.
(For more on the graph metric on trees, see~\cite{Bu}.)

Going back to our lamination picture,
we can associate to each arc
a non-negative real numbers, called its weight.
Such a weighted collection is a {\it measured lamination\/}.
Then, for any pair of vertices $i,j$,
we define $d_{i,j}$ as the sum of the weights of all the arcs
that separate $i$ from $j$.
This again gives a metric that satisfies
the ultrametric relation.
In the dual (tree-metric) picture,
this corresponds to assigning weights to edges,
and measuring distance between leaves
by summing weights along the path between them
rather than merely counting the edges.

For an extensive generalization of the foregoing picture,
see~\cite{FST}.

\section{A variant} \label{sec-variant}

An open problem concerns a variant
of Conway and Coxeter's definition,
in which the frieze recurrence is replaced by the recurrence
$$
\begin{array}{ccc}
  & A &   \\[1.8ex]
B & C & D \\[1.8ex]
  & E &  
\end{array} \ \ : \ \ E = (BD-C)/A
$$
and its sideways version
$$
\begin{array}{ccc}
  & A &   \\[1.8ex]
B & C & D \\[1.8ex]
  & E &  
\end{array} \ \ : \ \ D = (AE+C)/B \ \ .
$$
We can construct arrays 
that have the same sort of symmetries as frieze patterns
by starting with a suitable zig-zag of entries
(where successive downwards steps can go left, right, or straight)
and proceeding from left to right.
E.g., consider the partial table
$$
\begin{array}{ccccccc}
...& 1 & 1 & 1 & 1 & 1 & ...\\
   &   & A & D & x &   &    \\
   & B & E & y &   &   &    \\
   &   & C & F & z &   &    \\
...& 1 & 1 & 1 & 1 & 1 & ...
\end{array}
$$
where $A,...,F$ are pre-specified,
and where we compute
$y = (AC+E)/B$,
$x = (y+D)/A$,
$z = (y+F)/C$,
etc.
Then after exactly fourteen iterations of the procedure,
one gets back the original numbers (in their original order).
Moreover, along the way one sees
Laurent polynomials with positive coefficients.

Define a ``double zig-zag''
to be a subset of the entries of an $(n-2)$-rowed table
consisting of a pair of adjacent entries in each of the middle $n-4$ rows,
such that the pair in each row is displaced 
with respect to the pair in the preceding and succeeding rows
by at most one position.
(Thus the entries $A, B, C, D, E, F$ in the previous table 
form a double zig-zag,
as do the entries $D, E, F, x, y, z$.)

{\sc Conjecture}:
Given any assignment of formal weights
to the $2(n-4)$ entries in a double zig-zag
in an $(n-2)$-rowed table,
there is a unique assignment of rational functions
to all the entries in the table
so that the variant frieze relation is satisfied.
These rational functions of the original $2(n-4)$ variables
have glide-reflection symmetry that gives each row period $2n$.
Furthermore, each of the rational functions in the table
is a Laurent polynomial with positive coefficients.

There ought to be a way to prove this
by constructing the numerators of these Laurent polynomials
as sums of weights of perfect matchings of some suitable graph
(or perhaps sums of weights of combinatorial objects
more general than perfect matchings),
and the numerators undoubtedly contain abundant clues
as to how this can be done.

For $n=5,6,7,8$,
it appears that the number of positive integer arrays
satisfying the variant frieze relation
is 1, 5, 51, 868 (respectively).
This variant of the Catalan sequence
does not appear to have been studied before.
However, it should be said that these numbers 
were not computed in a rigorous fashion.
Indeed, it is conceivable that beyond some point,
the numbers becomes infinite
(i.e., for some $n$ there could be 
infinitely many $(n-2)$-rowed positive integer arrays
satisfying the variant frieze relation).

Dean Hickerson~\cite{Hi}~has shown 
that any $(n-2)$-rowed array that begins and ends with a row of 1's
and satisfies the variant frieze relation in between
has glide-reflection symmetry and period $2n$.
This implies that if one generates such a variant frieze pattern
starting with a double zig-zag of 1's,
one gets a periodic array of positive rational numbers.
However, it is not apparent that one can modify Hickerson's
(purely algebraic) proof to show that these rational numbers are integers.
Furthermore, if one uses formal weights instead of 1's,
Hickerson's argument does not seem to show
that the resulting rational functions are Laurent polynomials
(let alone that the Laurent polynomials have positive coefficients).


\section{Markoff numbers} \label{sec-markoff}

A {\it Markoff triple} is a triple $(x,y,z)$ of positive integers
satisfying $x^2+y^2+z^2=3xyz$, such as the triple (2,5,29).
A {\it Markoff number} is a positive integer
that occurs in at least one such triple.

Writing the Markoff equation as $z^2 - (3xy) z + (x^2+y^2) = 0$,
a quadratic equation in $z$, we see that if $(x,y,z)$ is a Markoff triple,
then so is $(x,y,z')$, where $z' = 3xy-z = (x^2+y^2)/z$,
the other root of the quadratic in $z$.
($z'$ is positive because $z'=(x^2+y^2)/z$,
and is an integer because $z'=3xy-z$.)
Likewise for $x$ and $y$.

The following claim is well-known
(for an elegant proof, see~\cite{Ba}):
Every Markoff triple $(x,y,z)$
can be obtained from the Markoff triple $(1,1,1)$
by a sequence of such exchange operations,
in fact, by a sequence of exchange operations
that leaves two numbers alone and increases the third.
E.g., $(1,1,1) \rightarrow (2,1,1) \rightarrow (2,5,1)$ $\rightarrow (2,5,29)$.

Create a graph whose vertices are the Markoff triples
and whose edges correspond to the exchange operations
$(x,y,z) \leftrightarrow (x',y,z)$,
$(x,y,z) \leftrightarrow (x,y',z)$,
$(x,y,z) \leftrightarrow (x,y,z')$
where
$x' = \frac{y^2+z^2}{x}$,
$y' = \frac{x^2+z^2}{y}$,
$z' = \frac{x^2+y^2}{z}$.
This 3-regular graph is connected (see the claim in the preceding paragraph), 
and it is not hard to show that it is acyclic.
Hence the graph is the 3-regular infinite tree.

This tree can be understood as the dual of the
triangulation of the upper half plane
by images of the modular domain under the action of the modular group.
Concretely, we can describe this picture by using Conway's terminology of
``lax vectors'', ``lax bases'', and ``lax superbases'' (\cite{C}).

A {\it primitive} vector $\vu$ in a lattice $L$
is one that cannot be written as $k\vv$
for some vector $\vv$ in $L$, with $k>1$.
A {\it lax vector} is a primitive vector defined only up to sign;
if $\vu$ is a primitive vector,
the associated lax vector is written $\pm \vu$.
A {\it lax base} for $L$ is a set of two lax vectors $\{\pm \vu, \pm \vv\}$
such that $\vu$ and $\vv$ form a basis for $L$.
A {\it lax superbase} for $L$ is a set of three lax vectors
$\{\pm \vu, \pm \vv, \pm \vw\}$ such that $\pm\vu\pm\vv\pm\vw=\vzero$
(with appropriate choice of signs)
and any two of $\vu,\vv,\vw$ form a basis for $L$. 

Each lax superbase
$\{\pm \vu, \pm \vv, \pm \vw\}$
contains the three lax bases
$\{\pm \vu, \pm \vv\}$, $\{\pm \vu, \pm \vw\}$, $\{\pm \vv, \pm \vw\}$
and no others.
In the other direction, each lax base
$\{\pm \vu, \pm \vv\}$
is in the two lax superbases
$\{\pm \vu, \pm \vv, \pm (\vu+\vv)\}, \ \{\pm \vu, \pm \vv, \pm (\vu-\vv)\}$
and no others.

The {\it topograph} is the graph whose vertices are lax superbases
and whose edges are lax bases, 
where each lax superbase is incident with the three lax bases in it.
This gives a 3-valent tree
whose vertices correspond to the lax superbases of $L$,
whose edges correspond to the lax bases of $L$,
and whose ``faces'' correspond to the lax vectors in $L$.

The lattice $L$ that we will want to use is the triangular lattice
$L = \{(x,y,z) \in \Z^3 \ : \ x+y+z=0\}$
(or $\Z^3 / \Z\vv$ where $\vv=(1,1,1)$, if you prefer).

Using this terminology, it is now possible to state 
the main idea of this section
(with details and proof to follow):
Unordered Markoff triples are associated with 
lax superbases of the triangular lattice,
and Markoff numbers are associated with 
lax vectors of the triangular lattice.
For example, the unordered Markoff triple $2,5,29$
corresponds to the lax superbase
$\{\pm \vu,\pm \vv,\pm \vw\}$
where $\vu = \vec{OA}$, $\vv = \vec{OB}$, and $\vw = \vec{OC}$,
with $O$, $A$, $B$, and $C$ forming a fundamental parallelogram
for the triangular lattice, as shown in Figure 19.
The Markoff numbers 1, 2, 5, and 29 correspond to the primitive vectors
$\vec{AB}$, $\vec{OA} = \vec{BC}$, 
$\vec{OB} = \vec{AC}$, and $\vec{OC}$.

\begin{center}
\begin{pspicture}*(0.5,0.4)(9.5,6.528) 
\psline(0.5,1.732)(9.5,1.732)
\psline(0.5,3.464)(9.5,3.464)
\psline(0.5,5.196)(9.5,5.196)
\psline(-1.5,0.866)(1.5,6.062) 
\psline(0.5,0.866)(3.5,6.062)
\psline(2.5,0.866)(5.5,6.062)
\psline(4.5,0.866)(7.5,6.062)
\psline(6.5,0.866)(9.5,6.062)
\psline(8.5,0.866)(11.5,6.062) 
\psline(-1.5,6.062)(1.5,0.866) 
\psline(0.5,6.062)(3.5,0.866)
\psline(2.5,6.062)(5.5,0.866)
\psline(4.5,6.062)(7.5,0.866)
\psline(6.5,6.062)(9.5,0.866)
\psline(8.5,6.062)(11.5,0.866) 
\rput(1,1.732){{\huge $O$}}
\rput(4,3.464){{\huge $A$}}
\rput(6,3.464){{\huge $B$}}
\rput(9,5.196){{\huge $C$}}
\end{pspicture}
\end{center}
\noindent

\begin{center}
Figure 19. A fundamental parallelogram.
\end{center}

To find the Markoff number associated with a primitive vector $\vec{OX}$,
take the union $R$ of all the triangles that segment $OX$ passes through.
The underlying lattice provides a triangulation of $R$.
E.g., for the vector $\vu = \vec{OC}$ from Figure 19,
the triangulation is as shown in Figure 20.
\begin{center}
\begin{pspicture}(0,0)(10,3.464) 
\psline(0,0.000)(4,0.000)
\psline(1,1.732)(7,1.732)
\psline(4,3.464)(8,3.464)
\psline(0,0.000)(1,1.732)
\psline(2,0.000)(4,3.464)
\psline(4,0.000)(6,3.464)
\psline(7,1.732)(8,3.464)
\psline(1,1.732)(2,0.000)
\psline(3,1.732)(4,0.000)
\psline(4,3.464)(5,1.732)
\psline(6,3.464)(7,1.732)
\rput(0,0.000){{\huge $O$}}
\rput(3,1.732){{\huge $A$}}
\rput(5,1.732){{\huge $B$}}
\rput(8,3.464){{\huge $C$}}
\psset{linestyle=dashed}
\psline(0,0)(8,3.464)
\end{pspicture}
\end{center}
\begin{center}
Figure 20. A Markoff snake.
\end{center}

Turn the triangulation into a planar bipartite graph as in section 2,
let $G(\vu)$ be the graph that results
from deleting vertices $O$ and $C$,
and let $M(\vu)$ be the number of perfect matchings of $G(\vu)$.
(If $\vu$ is a shortest vector in the lattice, put $M(\vu)=1$.)

\begin{theorem}[Carroll, Itsara, Le, Musiker, Price, and Viana
\cite{CILP} \cite{CP} \cite{I} \cite{M}]
If $\{\vu,\vv,\vw\}$ is a lax superbase of the triangular lattice,
then $(M(\vu),M(\vv),M(\vw))$ is a Markoff triple.
Every Markoff triple arises in this fashion.
In particular, 
if $\vu$ is a primitive vector, then $M(\vu)$ is a Markoff number,
and every Markoff number arises in this fashion.
\end{theorem}

\noindent
(The association of Markoff numbers with the topograph is not new;
what is new is the combinatorial interpretation
of the association, by way of perfect matchings.)

\begin{proof}
The base case, with $$(M(\ve_1),M(\ve_2),M(\ve_3)) = (1,1,1),$$ is clear.
The only non-trivial part of the proof is the verification that
$$M(\vu+\vv) = (M(\vu)^2 + M(\vv)^2) / M(\vu-\vv).$$ 
E.g., in Figure 21, we need to verify that
$$M(\vec{OC}) M(\vec{AB}) = M(\vec{OA})^2 + M(\vec{OB})^2.$$
But if we rewrite the desired equation as
$$M(\vec{OC}) M(\vec{AB}) = M(\vec{OA}) M(\vec{BC}) + M(\vec{OB}) M(\vec{AC})$$
we see that this is just Kuo's lemma
(see the proof of Theorem \ref{thm-cp}).
\begin{center}
\begin{pspicture}(0,0)(10,3.464) 
\psline(0,0.000)(4,0.000)
\psline(1,1.732)(7,1.732)
\psline(4,3.464)(8,3.464)
\psline(0,0.000)(1,1.732)
\psline(2,0.000)(4,3.464)
\psline(4,0.000)(6,3.464)
\psline(7,1.732)(8,3.464)
\psline(1,1.732)(2,0.000)
\psline(3,1.732)(4,0.000)
\psline(4,3.464)(5,1.732)
\psline(6,3.464)(7,1.732)
\rput(0,0.000){{\huge $O$}}
\rput(3,1.732){{\huge $A$}}
\rput(5,1.732){{\huge $B$}}
\rput(8,3.464){{\huge $C$}}
\end{pspicture}
\end{center}
\begin{center}
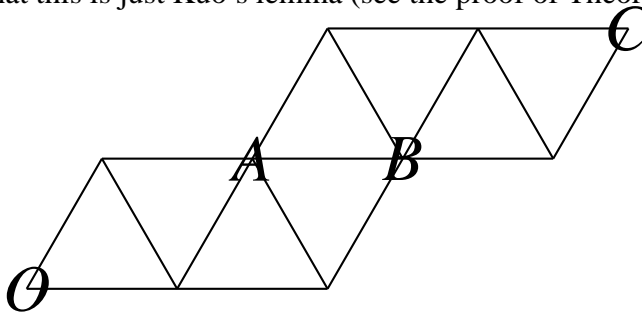
Figure 21. Kuo condensation for snakes.
\end{center}
\end{proof}


Remark: Some of the work done by Carroll et al.\ during the years of the 
Research Experiences in Algebraic Combinatorics at Harvard (2001 to 2003)
used the square lattice picture of section~\ref{sec-snake};
this way of interpreting the Markoff numbers combinatorially
was actually conjectured first, in 2001--2002, by Musiker,
and subsequently proved in 2002--2003 by Itsara, Le, Musiker, and Viana
(see~\cite{M},~\cite{I}, and~\cite{CILP},
and section~\ref{sec-snake} of this article).
Subsequently, the group's first combinatorial model
for frieze patterns, found by Price,
involved paths rather than perfect matchings.
It is reminiscent of, but apparently distinct from,
the paths model considered in section~\ref{sec-snake}.
Carroll turned Price's paths model into a perfect matchings model,
which made it possible to arrive at the snake-graph model
via a different route.

More generally, one can put $M(\ve_1)=x$, $M(\ve_2)=y$, and $M(\ve_3)=z$
(with $x,y,z>0$)
and recursively define $$M(\vu+\vv) = (M(\vu)^2 + M(\vv)^2) / M(\vu-\vv).$$
Then for all primitive vectors $\vu$,
$M(\vu)$ is a Laurent polynomial in $x,y,z$;
that is, it can be written in the form $P(x,y,z)/$ $x^a y^b z^c$,
where $P(x,y,z)$ is an ordinary polynomial
in $x,y,z$ (with non-zero constant term).
The numerator $P(x,y,z)$ of each Markoff polynomial is the sum of the weights
of all the perfect matchings of the graph $G(\vu)$, 
where edges have weight $x$, $y$, or $z$ according to orientation.
The triples $X=M(\vu)$, $Y=M(\vv)$, $Z=M(\vw)$
of rational functions associated with lax superbases
are solutions of the equation
$$X^2+Y^2+Z^2 = \frac{x^2+y^2+z^2}{xyz} XYZ.$$
Theorem~\ref{thm-frieze} implies that these numerators $P(x,y,z)$
are polynomials with positive coefficients.
This proves the following theorem:

\begin{theorem} \label{thm-triples}
Consider the initial triple $(x,y,z)$,
along with every triple of rational functions in $x$, $y$, and $z$
that can be obtained from the initial triple
by a sequence of operations of the form
$(X,Y,Z) \mapsto (X',Y,Z)$,
$(X,Y,Z) \mapsto (X,Y',Z)$, or
$(X,Y,Z) \mapsto (X,Y,Z')$,
where $X' = (Y^2+Z^2)/X$, $Y' = (X^2+Z^2)/Y$, and $Z' = (X^2+Y^2)/Z$.
Every rational function of $x$, $y$, and $z$ that occurs in such a triple
is a Laurent polynomial with positive coefficients.
\end{theorem}

Fomin and Zelevinsky proved in~\cite{FZL}
(Theorem 1.10)
that the rational functions $X(x,y,z),Y(x,y,z),Z(x,y,z)$ 
are Laurent polynomials,
but their methods did not prove positivity.
An alternative proof of positivity, based on topological ideas,
was given by Dylan Thurston~\cite{T}.

It can be shown that 
if $\vu$ is inside the cone generated by $+\ve_1$ and $-\ve_3$,
then $a < b > c$ and $(c+1) \ve_1 - (a+1) \ve_3 = \vu$.
(Likewise for the other sectors of $L$.)
This implies that all the ``Markoff polynomials'' $M(\vu)$ are distinct
(aside from the fact that $M(\vu)$ $=M(-\vu)$),
and thus $M(\vu)(x,y,z)$ $\neq M(\vv)(x,y,z)$ 
for all primitive vectors $\vu \neq \pm \vv$
as long as $(x,y,z)$ lies in a dense $G_\delta$ set of real triples.
This fact can be used to show~\cite{T}
that, for a generic choice of hyperbolic structure on the once-punctured torus,
no two simple geodesics have the same length.
(It should be mentioned that for the specific choice $x=y=z=1$,
the distinctness of the numbers $M(\vu)(x,y,z)$ as $\vu$ varies
is the famous, and still unproved, ``unicity conjecture'' for Markoff numbers.)


A slightly different point of view of Markoff numbers
focuses on triangles rather than lax superbases:
Say that points $A$, $B$, and $C$ in the equilateral triangular lattice
form a ``fundamental triangle'' if the area of triangle $ABC$
equals the area of the equilateral triangles of which the lattice is composed.
For example, the points $A$, $B$, and $C$ in Figure 21
are the vertices of a fundamental triangle.
(If four points form a fundamental parallelogram for the lattice,
then any three of the four points form a fundamental triangle.)
By Pick's theorem, a triangle is fundamental
if and only if it has no lattice points in its interior
and no lattice points on its boundary other than its three vertices.
Let $A$, $B$, and $C$ form a fundamental triangle.
Define the ``triangulation distance'' $d(x,y)$ between two vertices $x$ and $y$
as $M(\vu)$ where $\vu$ is the vector from $x$ to $y$.
Then the triangulation distances $d(A,B)$, $d(A,C)$, and $d(B,C)$
form a Markoff triple, and every Markoff triples arises in this way.

We conclude by mentioning a special sequence of Markoff numbers,
obtained by following the tree 
along those branches that give greatest numerical increase:
1, 1, 2, 5, 29, 433, 37666, ...
This sequence was considered by Dana Scott (see~\cite{G}),
and satisfies the recurrence
$f(n) = (f(n-1)^2+f(n-2)^2)/f(n-3)$.
Using the $A$ and $B$ matrices
from Section~\ref{sec-snake},
we see that we can alternately characterize the numbers
as the upper-left entries in the sequence of matrices
$$
\left( \begin{array}{rr}
1 & 1 \\ 1 & 0 \end{array} \right) , 
\left( \begin{array}{rr}
1 & 0 \\ 0 & 1 \end{array} \right) , 
\left( \begin{array}{rr}
2 & 1 \\ 1 & 1 \end{array} \right) , 
\left( \begin{array}{rr}
5 & 2 \\ 2 & 1 \end{array} \right) , 
\left( \begin{array}{rr}
29 & 12 \\ 12 & 5 \end{array} \right) , 
\left( \begin{array}{rr}
433 & 179 \\ 179 & 74 \end{array} \right) , ...
$$
satisfying the multiplicative recurrence relation
$$M(n) = M(n-1) M(n-3)^{-1} M(n-1)$$
(note that the Fibonacci numbers satisfy 
the additive version of this recurrence).
Andy Hone~\cite{Ho}
has shown that $\log \frac{\log f(n)}{n}$ 
approaches $\log \frac{1+\sqrt{5}}{2}$
as $n \rightarrow \infty$.


\section{Other directions for exploration} \label{sec-other}

\subsection{Non-integer frieze-patterns}

Given that the original geometric context of frieze patterns
gives rise to arrays containing numbers that are not integers,
it seems fairly natural to try to extend 
the Conway-Coxeter theory to this broader setting.
Enumerative questions would be a good place to start.
One might for instance try to count all the frieze patterns of order $n$
whose entries are either (positive) integers or half-integers,
and see if the enumerating sequence
is any sort of known analogue of the Catalan sequence.
Also, since many geometric frieze patterns
involve (irrational) algebraic numbers,
it might also be natural to enumerate
frieze patterns with entries in a given number ring
(though this might not be so very natural after all:
consider that, in its original geometric context,
positivity of the entries of the frieze pattern
is a consequence of their metric interpretation,
whereas for algebraic number rings
positivity is not a very robust notion
since it depends on the embedding of the ring in $\R$).


\subsection{Non-fundamental triangles}

Suppose $A$, $B$, and $C$ are points in the lattice
such that line segments $AB$, $AB$, and $BC$
contain no lattice-points other than their endpoints,
so that the triangulation distances $d(A,B)$, $d(A,C)$, $d(B,C)$
are well-defined.
We have seen that if triangle $ABC$ contains no lattice points
in its interior,
then these distances satisfy the Markoff equation.
Can anything be said if this condition does not hold?
For instance, in a lattice made of equilateral triangles of side-length 1,
consider an equilateral triangle $ABC$ of side-length $\sqrt{3}$
containing one interior point.
The triangulation distances are all equal to 2,
and $(2,2,2)$ is not a Markoff triple.
Nevertheless, perhaps there is a different algebraic equation
that this triple satisfies. 
More precisely, there may be an algebraic relation
satisfied by the triangulation distances 
$x=d(B',C')$, $y=d(A',C')$, $z=d(A',B')$
where $A'B'C'$ is any image of $ABC$
under the joint action of $SL_2(\Z)$ (change of lattice-base) 
and $\Z^2$ (translation).
Indeed, the whole numbers $x,y,z$ satisfy the condition
that there exist other whole numbers $x',y',z'$
(namely, the triangulation distances
from the interior point to $A'$, $B'$, and $C'$ respectively)
such that $x'^2+y^2+z^2 = 3x'yz$, $x^2+y'^2+z^2 = 3xy'z$,
and $x^2+y^2+z'^2 = 3xyz'$,
and perhaps some sort of quantifier elimination procedure
would permit us to write this as a condition on just $x$, $y$, and $z$.
More broadly, perhaps each orbit of triangles
under the action of $SL_2(\Z)$ and $\Z^2$
gives rise to triples satisfying
a particular algebraic condition specific to that orbit.

\subsection{Other ternary cubics}

Neil Herriot (another member of REACH) 
showed~\cite{He} that if we replace the triangular lattice used above
by the tiling of the plane by isosceles right triangles
(generated from one such triangle by repeated reflection in the sides),
fundamental triangles give rise to triples $x,y,z$ of positive integers
satisfying either $$x^2+y^2+2z^2=4xyz$$ or $$2x^2+2y^2+z^2=4xyz.$$
(Note that these two Diophantine equations are essentially equivalent,
as the map $(x,y,z) \mapsto (x,y,2z)$ gives a bijection between
solutions to the former and solutions to the latter.)
For instance, if for any two vertices $X,Y$
we define the triangulation distance $d(X,Y)$ 
in analogy with the definition used before
(now using the isosceles right triangle lattice
in place of the equilateral triangle lattice),
then the points $O,A,B,C$ shown in Figure 22 satisfy
$d(A,B) = 1$, $d(O,A) = 1$, $d(B,C) = 2$,
$d(O,B) = 3$, $d(A,C) = 3$, and $d(O,C) = 11$,
corresponding to the solution triples
$(11)^2 + (3)^2 + 2(1)^2 = 4(11)(3)(1)$
and
$2(11)^2 + 2(3)^2 + (2)^2 = 4(11)(3)(2)$.
\begin{center}
\begin{pspicture}(1,1)(9,7)
\psline(2,2)(8,2)
\psline(2,4)(8,4)
\psline(2,6)(8,6)
\psline(2,2)(2,6)
\psline(4,2)(4,6)
\psline(6,2)(6,6)
\psline(8,2)(8,6)
\psline(2,2)(6,6)
\psline(6,2)(8,4)
\psline(6,2)(2,6)
\psline(8,4)(6,6)
\rput(1.8,2.2){{\Large $O$}}
\rput(3.8,4.2){{\Large $A$}}
\rput(5.8,4.2){{\Large $B$}}
\rput(7.8,6.2){{\Large $C$}}
\end{pspicture}
\end{center}
\begin{center}
Figure 22.  Herriot's theorem.
\end{center}
More specifically, Herriot showed that
if $ABC$ is a fundamental triangle,
then the triangulation distances
$d(A,B)$, $d(A,C)$, $d(B,C)$ satisfy
$$d(A,B)^2 + 2d(A,C)^2 + 2d(B,C)^2 = 4d(A,B)d(A,C)d(B,C)$$
or
$$2d(A,B)^2 + d(A,C)^2 + d(B,C)^2 = 4d(A,B)d(A,C)d(B,C)$$
according to whether the vertices $A,B,C$
have respective degrees 4,4,8 or 8,8,4.
(One can check that a fundamental triangle
cannot have all three vertices of degree 4
or all three vertices of degree 8.)
A related observation is that
if $OACB$ is a fundamental parallelogram
with $O$ and $A$ of degree 8 and $B$ and $C$ of degree 4,
then $d(B,C) = 2d(O,A)$.

Herriot's result, 
considered in conjunction with the result on Markoff numbers,
raises the question of whether
there might be some more general combinatorial approach
to ternary cubic equations of similar shape.

Rosenberger~\cite{R}
showed that there are exactly three ternary cubic equations
of the shape $ax^2 + by^2 + cz^2 = (a+b+c)xyz$
for which all the positive integer solutions
can be derived from the solution $(x,y,z)=(1,1,1)$ by means
of the exchange operations 
$(x,y,z) \rightarrow (x',y,z)$,
$(x,y,z) \rightarrow (x,y',z)$, and
$(x,y,z) \rightarrow (x,y,z')$, with 
$x'=(by^2+cz^2)/ax$,
$y'=(ax^2+cz^2)/by$, and
$z'=(ax^2+by^2)/cz$.
These three ternary cubic equations are
$$x^2+y^2+z^2=3xyz,$$
$$x^2+y^2+2z^2=4xyz,$$ and
$$x^2+2y^2+3z^2=6xyz.$$

Note that the triples of coefficients that occur here ---
(1,1,1), (1,1,2), and (1,2,3) ---
are precisely the triples that occur in the classification
of finite reflection groups in the plane.
Specifically, the ratios 1:1:1, 1:1:2, and 1:2:3
describe the angles of the three triangles ---
the 60-60-60 triangle, the 45-45-90 triangle, and the 30-60-90 triangle ---
that arise as the fundamental domains
of the three irreducible two-dimensional reflection groups.

Since the solutions to the ternary cubic $x^2+y^2+z^2=3xyz$
describe properties of the tiling of the plane by 60-60-60 triangles,
and solutions to the ternary cubic $x^2+y^2+2z^2=4xyz$
describe properties of the tiling of the plane by 45-45-90 triangles,
the solutions to the ternary cubic $x^2+2y^2+3z^2=6xyz$
``ought'' to be associated with some combinatorial model 
involving the reflection-tiling of the plane by 30-60-90 triangles.
Unfortunately, the most obvious approach
(based on analogy with the 60-60-60 and 45-45-90 cases) does not work.
So we are left with two problems that may or may not be related:
first, to find a combinatorial interpretation
for the integers (or, more generally, the Laurent polynomials)
that arise from solving the ternary cubic $x^2+2y^2+3z^2=6xyz$;
and second, to find algebraic recurrences
that govern the integers (or, more generally, the Laurent polynomials)
that arise from counting (or summing the weights of)
perfect matchings of graphs
derived from the reflection-tiling of the plane by 30-60-90 triangles.

If there is a way to make the analogy work,
one might seek to extend the analysis to other ternary cubics.
It is clear how this might generalize on the algebraic side. 
On the geometric side, one might drop the requirement
that the triangle tile the plane by reflection,
and insist only that each angle be a rational multiple of 360 degrees.
There is a relatively well-developed theory
of ``billiards flow'' in such a triangle
(see e.g.~\cite{KS})
where a particle inside the triangle
bounces off the sides following the law of reflection
(angle of incidence equals angle of reflection)
and travels along a straight line in between bounces.
The path of such a particle can be unfolded
by repeatedly reflecting the triangular domain
in the side that the particle is bouncing off of,
so that the unfolded path of the particle
is just a straight line in the plane.
Of special interest in the theory of billiards
are trajectories joining a corner to a corner
(possibly the same corner or possibly a different one);
these are called saddle connections.
The reflected images of the triangular domain 
form a triangulated polygon,
and the saddle connection is a combinatorial diagonal of this polygon.
It is unclear whether the combinatorics of such triangulations
might contain dynamical information about the billiards flow,
but if this prospect were to be explored,
enumeration of matchings on the derived bipartite graphs
would be one thing to try.

\subsection{More variables}

Another natural variant of the Markoff equation $x^2+y^2+z^2=3xyz$
is the equation $w^2+x^2+y^2+z^2=4wxyz$
(one special representative of a broader class
called Markoff-Hurwitz equations; see~\cite{Ba}).
The Laurent phenomenon applies here too:
the four natural exchange operations
convert an initial formal solution $(w,x,y,z)$
into a quadruple of Laurent polynomials.
(This is a special case of Theorem 1.10 in~\cite{FZL}.)

Furthermore, the coefficients of these Laurent polynomials
appear to be positive, although this has not been proved.

The numerators of these Laurent polynomials ought to be weight-enumerators
for some combinatorial model, but it is unclear
how to reverse-engineer the combinatorial model from the Laurent polynomials.


\bibliographystyle{amsalpha}

\end{document}